\documentclass[reqno]{amsart}
\usepackage{amssymb}
\usepackage{graphicx}

\usepackage[usenames, dvipsnames]{color}
\usepackage{verbatim}
\usepackage{mathrsfs}
\usepackage{bm}
\usepackage{cite}
\usepackage{bbm}

\numberwithin{equation}{section}

\newtheorem{theorem}{Theorem}[section]

\newtheorem{lemma}[theorem]{Lemma}
\newtheorem{prop}[theorem]{Proposition}

\theoremstyle{definition}
\newtheorem{remark}[theorem]{Remark}

\theoremstyle{definition}

\theoremstyle{definition}

\makeatletter
\def\dashint{\operatorname%
{\,\,\text{\bf-}\kern-.98em\DOTSI\intop\ilimits@\!\!}}
\makeatother

\def\\det{\text{\det}}

\def\Xint#1{\mathchoice
 {\XXint\displaystyle\textstyle{#1}}%
 {\XXint\textstyle\scriptstyle{#1}}%
 {\XXint\scriptstyle\scriptscriptstyle{#1}}%
 {\XXint\scriptscriptstyle\scriptscriptstyle{#1}}%
 \!\int}
\def\XXint#1#2#3{{\setbox0=\hbox{$#1{#2#3}{\int}$}
  \vcenter{\hbox{$#2#3$}}\kern-.5\wd0}}

\def\dashint{\Xint-}

\def\.5{\frac{1}{2}}

\newcommand{\RN}[1]{%
  \textup{\uppercase\expandafter{\romannumeral#1}}%
}

\renewcommand{\epsilon}{\varepsilon}

\newcounter{marnote}

\begin{document}

\title[Higher integrability for parabolic $p$-Laplace equations]{Higher integrability for parabolic $p$-Laplace equations with gradient-dependent sources}

\author[Y. Li]{Yan Li}
\address[Y. Li] {School of Mathematical Sciences, Beijing Normal University, Beijing 100875, China. }
\email{yanli25@mail.bnu.edu.cn}

\author[Z.W. Zhao]{Zhiwen Zhao}

\address[Z.W. Zhao]{School of Mathematics and Physics, University of Science and Technology Beijing, Beijing 100083, China.}
\email{zwzhao365@163.com;\,zwzhao365@gmail.com}

%\footnote{}

\date{\today} % delete this line to display the current date

%%% BEGIN DOCUMENT

\maketitle
%\tableofcontents
\begin{abstract}
This paper aims to establish higher gradient integrability away from the initial time for solutions to parabolic $p$-Laplace equations with general gradient-dependent lower-order terms throughout the singular and degenerate regimes.\ For this purpose, we develop an elliptic reconstruction argument based on an $L^2$-valued temporal H\"{o}lder estimate, thereby transferring the problem to elliptic estimates and quantitative control of the reconstruction error.\ In particular, our approach yields an explicit higher-integrability exponent.\ We further construct counterexamples illustrating that bounded continuous initial data alone are insufficient to ensure higher gradient integrability up to the initial time. Finally, we provide counterexamples demonstrating that gradient-dependent lower-order terms can prevent temporal Lipschitz regularity.

\end{abstract}

\maketitle
%\date{}
%\maketitle
%{\bf Abstract}

%{\bf{Keywords:}} H\"{o}lder estimates; singular parabolic $p$-Laplace equations; temporal Lipschitz regularity
%
%\noindent{\bf{MSC numbers}}: {35K92; 35B51; 35B65.}

\section{Introduction}
Suppose that $\Omega\subset\mathbb R^n$ is a smooth bounded domain with $n\geq2$. For any $T>0$, set $\Omega_T=\Omega\times(-T,0]$. Define
\begin{equation}\label{eq:spatial-target}
	q=\frac{np}{n-1},\quad\text{for } p\in(1,2)\cup(2,\infty),\;n\geq2,
\end{equation}
and
\begin{equation*}
A_{p}(\xi)=|\xi|^{p-2}\xi,\quad\text{for }\xi\in\mathbb{R}^{n},\;p\in(1,2)\cup(2,\infty).
\end{equation*}
Assume that
 \begin{align}\label{Tr01}
 g\in L^\infty(\partial \Omega)\cap W^{1-1/q,q}(\partial \Omega).
 \end{align}
This paper is concerned with the higher integrability of the spatial gradients of weak solutions to the following initial-boundary value problem for a parabolic $p$-Laplace equation with a gradient-dependent lower-order term:
\begin{equation}\label{PO001}
	\begin{cases}
		\partial_{t}u-\mathrm{div}A_{p}(\nabla u)=b(x,t,u,\nabla u),&\text{in }\Omega_{T},\\
		u=g(x)\;\,\text{or}\;\,
		A_{p}(\nabla u)\cdot\nu=0,&\text{on }\partial \Omega\times(-T,0],\\
		u(x,-T)=\phi(x),&\text{in }\Omega,
	\end{cases}
\end{equation}
where $\phi\in L^{\infty}(\Omega)$, $\nu$ is the outward unit normal to $\partial\Omega$. The two boundary conditions are imposed separately. The measurable lower-order term $b:\Omega_{T}\times\mathbb{R}^{n+1}\rightarrow\mathbb{R}$ satisfies
\begin{equation}\label{E02}
	|b(x,t,u,\nabla u)|\leq L|\nabla u|^{p-1}+\psi(x,t),\quad\psi\geq0,\quad \psi\in L^{\frac{pl_{0}}{p-1}}(\Omega_{T}),
\end{equation}
where $L>0$ and $l_{0}>\frac{n+p}{p}$. A function $u\in C([-T,0];L^2(\Omega)) \cap L^p((-T,0);W^{1,p}(\Omega))$ is called a weak solution of \eqref{PO001} whenever the following variational identity holds: for every $-T\leq t_{1}<t_{2}\leq0$,
\begin{align*}
	&\int_{\Omega}u\varphi\,dx\Big|_{t_{1}}^{t_{2}}+\int_{t_{1}}^{t_{2}}\int_{\Omega}\left(-u\partial_{t}\varphi
	+A_{p}(\nabla u)\cdot\nabla\varphi\right)dxdt=\int_{t_{1}}^{t_{2}}\int_{\Omega}b(x,t,u,\nabla u)\varphi dxdt,
\end{align*}
for all $\varphi\in W^{1,2}((-T,0);L^{2}(\Omega))\cap L^{p}((-T,0);W^{1,p}_{0}(\Omega))$ in the Dirichlet case and $\varphi\in W^{1,2}((-T,0);L^{2}(\Omega))\cap L^{p}((-T,0);W^{1,p}(\Omega))$ in the homogeneous Neumann case. For $p>2$, the parabolic $p$-Laplace equation
is said to be degenerate, since its modulus of ellipticity $|\nabla u|^{p-2}$ vanishes wherever $\nabla u=0$. In contrast, for $1<p<2$, the equation is termed singular because the same modulus becomes unbounded as $|\nabla u|\to0$.

The self-improvement of gradient integrability is a basic phenomenon in both elliptic and parabolic regularity theory.\ For nonlinear elliptic systems, a classical higher-integrability result was first established by Meyers and Elcrat \cite{ME1975}; see also \cite{G1983,S1980}. The standard argument combines Caccioppoli and Sobolev--Poincar\'e inequalities to obtain a reverse H\"older estimate for the gradient.\ A suitable version of Gehring's lemma then provides the crucial self-improvement, yielding integrability beyond the natural energy exponent.\ This self-improvement principle was originally developed by Gehring \cite{G1973} in his proof of higher integrability for the derivatives of quasiconformal mappings.\ Analogous results for parabolic systems with quadratic growth ($p=2$) were obtained by Giaquinta and Struwe \cite{GS1982}; see also \cite{C1993}. 

However, for nonlinear parabolic systems of $p$-Laplacian type with $p\neq2$, the mismatch between the quadratic time-energy term and the spatial dissipation of $p$-growth obstructs a direct extension of the elliptic argument on standard parabolic cylinders. A fundamental breakthrough was achieved by Kinnunen and Lewis \cite{KL2000}, who addressed this scaling mismatch through intrinsic cylinders adapted to the local energy scale.\ On these cylinders, they established suitably normalized reverse H\"older estimates, which, together with a covering argument and Gehring-type self-improvement, yield local gradient higher integrability in the supercritical range $p>\frac{2n}{n+2}$.\ The concept of intrinsic cylinders can be traced back to the work of DiBenedetto and Friedman \cite{DF198501} on the H\"{o}lder continuity of solutions (see also \cite{D1993}). Under stronger a priori integrability assumptions on the solution, notably $u\in L^r_{\mathrm{loc}}$ with $r>2$ and $n(p-2)+pr>0$, B\"{o}gelein, Duzaar, Gianazza and Liao \cite{BDGL2026} established local boundedness and gradient higher integrability in the critical and subcritical ranges.\ Their argument combines quantitative local boundedness with mixed intrinsic scaling adapted to both the solution and its gradient to derive the required reverse H\"{o}lder estimates.\ With regard to further studies on higher gradient integrability for other nonlinear parabolic equations, as well as extensions to very weak solutions, we refer to \cite{BDGL2025,KKM2023,BDKS2020,GS2019,KL2002} and the references therein.

Unlike the above approach based on reverse H\"{o}lder inequalities, in this paper we develop a new route to higher integrability that starts from the establishment of an $L^{2}$-valued temporal H\"older estimate.\ We exploit the resulting temporal regularity through elliptic reconstruction on time intervals, reducing the parabolic higher-integrability problem to elliptic estimates and quantitative control of the reconstruction error.\ Temporal H\"older continuity controls this error, while elliptic regularity controls the reconstructed gradients.\ A combination of these bounds leads to an estimate for gradient superlevel sets, with the time scale chosen to balance the two contributions. In particular, this quantitative reduction allows us to obtain an explicit decay estimate for the distribution function of the gradient and thus establishes higher integrability with an explicit exponent.

%Recently, Hao, Li, and Zhao \cite{HLZ202602} obtained space-time $W^{1,\infty}_{x,t}$ estimates away from the initial time for the classical parabolic $p$-Laplace equations with bounded nonnegative time-independent forcing and initial data. In particular, their estimates cover both the singular and degenerate ranges with explicit dependence on the prescribed data.

Define
\begin{align}\label{E07}
	\alpha=
	\begin{cases}
		\frac{1}{p+1},&1<p<2,\\
		\frac{1}{2p-1},&p>2,	
	\end{cases}	\quad \beta=	
	\begin{cases}
		\frac{1}{p},&1<p<2,\\
		\frac{1}{2(p+1)},&p>2,	
	\end{cases}	
\end{align}	
and
\begin{align}
	\sigma=&\alpha\min\{1,2/p\},\quad
	\gamma=\min\bigg\{1,\frac2n+\frac{2(p-1)}q\bigg\},\quad\kappa=q-p=\frac p{n-1},
	\label{eq:time-exponents}\\
	\vartheta=&\frac{1-\alpha\gamma}{p-1},\quad
	\mu=\max\{\kappa/p,q\vartheta\},\quad
	\varepsilon_*=\frac{\sigma\kappa}{\sigma+\mu}.
	\label{eq:gain-exponents}
\end{align}

Throughout the paper, $C$ represents a positive constant that may vary from line to line and depends only on the prescribed data, including $n$, $p$, $L$,
$\|\phi\|_{L^{\infty}(\Omega)}$, $\|g\|_{L^{\infty}(\partial\Omega)}$, $\|g\|_{W^{1-1/q,q}(\partial\Omega)}$, and $\|\psi\|_{L^{pl_{0}/(p-1)}(\Omega_T)}$.

We now state our main results.
\begin{theorem}\label{thm:main}
Fix $p\in(1,2)\cup(2,\infty)$ and $n\geq2$. Let $u$ be a weak solution of \eqref{PO001} with $\Omega_{T}=B_{1}\times(-1,0]$. Then one can find a constant $L_{\ast}\in(0,1)$, depending solely on $n$ and $p$, such that if $0<L\leq L_{\ast}$, for every $0<\varepsilon<\varepsilon_*$ and $0<\delta\ll1$, 
	\begin{equation}\label{eq:strong-conclusion}
	\int_{B_{1}\times(-1+\delta,0)}|\nabla u|^{p+\varepsilon}dxdt\leq C(\varepsilon_{\ast}-\varepsilon)^{-1}\delta^{-\beta\gamma q/(p-1)},
	\end{equation}
where $\beta$, $\gamma$ and $\varepsilon_{\ast}$ are given by \eqref{E07}, \eqref{eq:time-exponents} and \eqref{eq:gain-exponents}, respectively.
\end{theorem}

\begin{remark}
By picking $\varepsilon=\varepsilon_{\ast}/2$, we obtain the following explicit higher integrability exponent:
	\begin{align}\label{eq:explicit-epsilon}
	 p_{\sharp}=p+\varepsilon,\quad\varepsilon=\frac{\sigma p}{2\left[\sigma(n-1)+np(1-\alpha\gamma)/(p-1)\right]}.
	\end{align}	
For $n\geq2$, denote $p_{\mathrm{c}}=2-2/n.$ The gain in \eqref{eq:explicit-epsilon} becomes
\begin{align*}
\varepsilon(p,n)=
\begin{cases}
	\displaystyle\frac{p(p-1)}{2(np^2-p+n-1)},
	&1<p<p_{\mathrm{c}},\;n>2,\\[2.6mm]
	\displaystyle\frac{p(p-1)}
	{2[np^2+(n-1)(p-1)]},
	&p_{\mathrm{c}}\le p<2,\\[2.6mm]
	\displaystyle\frac{p}{2(np^2+n-1)},
	&p>2.
\end{cases}
\end{align*}
Direct differentiation shows that for fixed $n$, the chosen gain is strictly increasing on $(1,2)$ and strictly decreasing on $(2,\infty)$; for fixed $p$, it strictly decreases with increasing dimension $n$. Consequently, we have
\begin{align*}
\sup_{\substack{p>1\\p\neq2}}\varepsilon(p,n)=\lim_{p\to2}\varepsilon(p,n)=\frac{1}{5n-1}.
\end{align*}
Thus lower dimensions and values of $p$ close to $2$ give larger gains.\ At the two extremes, we obtain that for each fixed $n\geq2$, 
\begin{align*}
\varepsilon(p,n)=
\begin{cases}	
\frac{p-1}{4(n-1)}+O\bigl((p-1)^2\bigr),
	&\text{as }p\rightarrow1,\\[0.5mm]
\frac1{2np}+O(p^{-3}),&\text{as }p\to\infty,
\end{cases}
\end{align*}
where the implicit constants may depend on $n$. This implies that the upper bound in \eqref{eq:strong-conclusion} blows up as $p\rightarrow 1$ or $p\rightarrow\infty.$
	
\end{remark}	

\begin{remark}
A direct computation reads that
\begin{equation*}
	\frac{\beta\gamma q}{p-1}=
	\begin{cases}
		\displaystyle\frac2p+\frac{2}{(n-1)(p-1)},
		&1<p<p_{\mathrm{c}},\;n>2,\\[2.6mm]
		\displaystyle\frac{n}{(n-1)(p-1)},
		&p_{\mathrm{c}}\le p<2,\\[2.6mm]
		\displaystyle\frac{np}{2(n-1)(p^2-1)},
		&p>2,
	\end{cases}
\end{equation*}
where $p_{\mathrm{c}}=2-2/n.$ For fixed $n$, this quantity strictly decreases on each side of $2$ with a downward jump:
\begin{align*}
\lim_{p\to2^-}\frac{\beta\gamma q}{p-1}=\frac{n}{n-1},
\quad
\lim_{p\to2^+}\frac{\beta\gamma q}{p-1}=\frac{n}{3(n-1)}.
\end{align*}
For fixed $p$, this quantity strictly decreases as $n$ increases. Moreover, we have
\begin{align*}	
	\frac{\beta\gamma q}{p-1}=
\begin{cases}	
	\frac{2}{(n-1)(p-1)}+O(1),
	&\text{as }p\rightarrow1,\\
\frac{n}{2p(n-1)}+O(p^{-3}),
	&\text{as }p\to\infty.
\end{cases}	
\end{align*}
Here the implied constants may depend on $n$.\ This shows that the factor $\delta^{-\beta\gamma q/(p-1)}$ in the upper bound of \eqref{eq:strong-conclusion} diverges to infinity as $p\to1$.

\end{remark}	

\begin{remark}
The result in \eqref{eq:strong-conclusion} is global in space and hold up to the terminal time, while the initial time remains separated by $\delta>0$. The fixed trace assumption in \eqref{Tr01} is essential for this global argument and bounded Dirichlet data alone do not suffice. In fact, if $u\in L^r((-1+\delta,0);W^{1,r}(B_1))$ for some $r>p$, it follows from the trace inequality (see p.\ 315 in \cite{B2011}) that
	\begin{align*}
	(1-\delta)\|g\|_{W^{1-1/r,r}(\partial B_1)}^{r}\leq& \int_{-1+\delta}^{0}\|\mathrm{Tr}u(\cdot,t)\|_{W^{1-1/r,r}(\partial B_{1})}^{p}dt\\
	\leq& C(n,r)\int_{-1+\delta}^{0}\|u(\cdot,t)\|_{W^{1,r}(B_1)}^{r}dt,
	\end{align*}
where $\mathrm{Tr}:W^{1,r}(B_{1})\longrightarrow W^{1-1/r,r}(\partial B_{1})$ is the trace operator. Here $W^{1-1/r,r}(\partial B_1)$ denotes the fractional Sobolev trace space.\ Hence some higher trace regularity is necessary for a global gain. 
\end{remark}

\begin{remark}
For problem \eqref{PO001} in $B_1\times(-1,0]$, the energy regularity and fixed Dirichlet condition already imply $g\in W^{1-1/p,p}(\partial B_1)$ by using the trace theorem. Without the additional trace regularity in \eqref{Tr01}, and under the remaining hypotheses, the same argument, together with the interior elliptic estimate of Lemma \ref{interior-elliptic-CZ} in the appendix, yields the following local higher-integrability estimate:
\begin{align*}
\|\nabla u\|^{p+\varepsilon}_{L^{p+\varepsilon}(B_{R}\times(-1+\delta,0))}\leq C_R (\varepsilon_{\ast}-\varepsilon)^{-1}\delta^{-\beta\gamma q/(p-1)},\quad\text{for any }0<R<1,
\end{align*}
with the same admissible $\varepsilon$. Here $C_R\rightarrow\infty$, as $R\rightarrow1$.

\end{remark}

\begin{remark}
	The conclusion of \eqref{eq:strong-conclusion} extends to $\Omega\times(-T,0]$, where $\Omega\subset\mathbb R^n$ is a bounded connected domain with smooth boundary and $T>0$. Applying the proof of Theorem \ref{thm:main} with minor modifications, we obtain that for every $0<\delta\ll T$ and $0<\varepsilon<\varepsilon_\ast$,
	\begin{align*}
	\|\nabla u\|^{p+\varepsilon}_{L^{p+\varepsilon}
		(\Omega\times(-T+\delta,0))}\leq C_{\Omega,T}(\varepsilon_\ast-\varepsilon)^{-1}\delta^{-\beta\gamma q/(p-1)}.
	\end{align*}
	The exponents are the same as in the unit-cylinder case. Besides the original data dependence, the constant $C_{\Omega,T}$ depends on $\Omega$ and $T$, but is independent of $\delta$ and $\varepsilon$. In addition, when the lower-order term is replaced by a bounded nonnegative time-independent forcing term, space-time $W^{1,\infty}_{x,t}$ estimates after any fixed positive waiting time were established in Theorem 1.5 of \cite{HLZ202602}.

\end{remark}

The proof of Theorem \ref{thm:main} relies crucially on the temporal H\"{o}lder estimate stated below.
\begin{theorem}\label{thm10}
For $p\in(1,2)\cup(2,\infty)$ and $n\geq2$, let $u$ denote the solution of \eqref{PO001} in $B_1\times(-1,0]$.\ Given an arbitrary $0<\delta\ll1$, there exists a constant $L_{\ast}=L_{\ast}(n,p)\in(0,1)$ such that if $0<L\leq L_{\ast}$, for any $-1+\delta\leq s<t\leq 0$,
\begin{align*}
	\|u(\cdot,t)-u(\cdot,s)\|_{L^{2}(B_{1})}\leq C\delta^{-\beta}|t-s|^{\alpha},
\end{align*}
where $\alpha$ and $\beta$ are given by \eqref{E07}.

\end{theorem}
\begin{remark}
When the lower-order term reduces to a bounded nonnegative time-independent forcing term, $L^\infty$-valued temporal Lipschitz estimates have been established in Theorem 4.1 of \cite{HLZ202601} for $1<p<2$ and Theorem 3.1 of \cite{HLZ202602} for $p>2$. When the lower-order term is instead a gradient-independent forcing term satisfying finite integrability conditions, $L^\infty$-valued temporal H\"{o}lder estimates have been obtained in Theorem 1.3 of \cite{LZ202601}, covering the supercritical singular range and the whole degenerate regime.	
\end{remark}

To show that bounded continuous initial data alone do not guarantee higher integrability of the spatial gradient up to the initial time, we consider the following two alternative problems:
\begin{equation}\label{eq:problem}
	\begin{cases}
		\partial_tu-\mathrm{div} A_p(\nabla u)=0,
		&\text{in }B_1\times(-1,0],\\
		u=0\;\,\text{or}\;\, A_p(\nabla u)\cdot\nu=0,
		&\text{on }\partial B_1\times(-1,0],\\
		u(\cdot,-1)=\phi,
		&\text{in }B_1,
	\end{cases}
\end{equation}
where $\nu$ denotes the outward unit normal. Let the symbol $\Subset$ stand for compact containment.
\begin{prop}\label{prop:counterexample}
Assume that $p>1$ and $n\geq2$.\ There exist a merely continuous function $\phi$ compactly supported in $B_{1/2}$ with $0\leq\phi\leq2$, and a cube $D\Subset B_{1/2}$ such that the energy solution of \eqref{eq:problem}, for either choice of boundary condition, satisfies
	\begin{equation*}
		0\leq u\leq2,\quad\int_{-1}^{0}\int_{B_1}|\nabla u|^{p}dxdt\leq\frac{1}{2}\|\phi\|_{L^2(B_1)}^{2},
	\end{equation*}
	but
	\begin{equation}\label{AG09}
		\int_{-1}^{-1+\tau}\int_{D}|\nabla u|^{p+\varepsilon}dxdt=\infty,\quad
		\text{for every }\varepsilon>0\text{ and }0<\tau\leq1.
	\end{equation}
	In particular,
	\begin{equation*}
		u\notin L^{p+\varepsilon}\bigl((-1,0);W^{1,p+\varepsilon}(B_1)\bigr),
		\quad\text{for any }\varepsilon>0.
	\end{equation*}
	
\end{prop}

\begin{remark}
The failure in \eqref{AG09} occurs in a fixed interior cube and on every time interval adjacent to the initial time.\ It therefore persists even with zero forcing, smooth homogeneous lateral boundary data, and a merely continuous, compactly supported initial datum.\ Remarkably, the same phenomenon occurs
already for $p=2$, namely, for the classical heat equation.\ Thus, the obstruction arises from the limited regularity of the initial datum rather than from the nonlinearity of the equation.	

\end{remark}

We next construct counterexamples to $L^2$-valued temporal Lipschitz endpoint for the equations with gradient-dependent
lower-order terms. To this end, we consider bounded energy solutions of
\begin{equation}\label{eq:model}
	\begin{cases}
		\partial_tu-\mathrm{div}A_p(\nabla u)
		=b(x,t,u,\nabla u),
		&\text{in }B_1\times(-1,0],\\
		u=g\;\,\text{or}\;\, A_p(\nabla u)\cdot\nu=0,
		&\text{on }\partial B_1\times(-1,0],\\
		u(\cdot,-1)=\phi,
		&\text{in }B_1,
	\end{cases}
\end{equation}
where $g$ is time-independent and
\begin{align*}
|b(x,t,u,\nabla u)|\leq L|\nabla u|^{p-1},\quad L>0.
\end{align*}
For any interval $I\subset(-1,0]$ of positive length, define
\begin{equation}\label{eq:lip}
	[u]_{\operatorname{Lip}(I;L^2(B_1))}
	:=
	\sup_{\substack{s,t\in I\\s<t}}
	\frac{\|u(\cdot,t)-u(\cdot,s)\|_{L^2(B_1)}}{t-s}
	\in[0,\infty].
\end{equation}
Temporal Lipschitz continuity with values in $L^2$ is equivalent to finiteness of this seminorm.

\begin{prop}[Dirichlet counterexample]\label{prop:D}
Set $p\in(1,2)\cup(2,\infty)$ and $n\geq2$.	For every $L>0$, there exist $c_{\scriptscriptstyle D}\in C(\overline{B}_1\times[-1,0])$ with $\|c_{\scriptscriptstyle D}\|_{L^\infty(B_1\times(-1,0])}\leq L$, and a positive bounded energy solution $v_{\scriptscriptstyle D}$ of the Dirichlet problem in \eqref{eq:model} with $b=c_{\scriptscriptstyle D}|\nabla v_{\scriptscriptstyle D}|^{p-1}$ and $g(x)=\phi(x)=2+x_1$.\ Moreover, $\frac12\le|\nabla v_{\scriptscriptstyle D}|\leq\frac32$ in $B_1\times(-1,0]$, and the source $c_{\scriptscriptstyle D}|\nabla v_{\scriptscriptstyle D}|^{p-1}$ is continuous and bounded. Nevertheless, $[v_{\scriptscriptstyle D}]_{\operatorname{Lip}(I;L^2(B_1))}=\infty$ for every $I=(-\frac12,-\frac12+\rho)$, $0<\rho\le\frac12$.
\end{prop}

\begin{prop}[Neumann counterexample]\label{prop:N}
Let $p\in(1,2)\cup(2,\infty)$ and $n\geq2$.	For every $L>0$, there exist $c_{\scriptscriptstyle N}\in C(\overline{B}_1\times[-1,0])$ with $\|c_{\scriptscriptstyle N}\|_{L^\infty(B_1\times(-1,0])}\leq L$, and a positive bounded energy solution $v_{\scriptscriptstyle N}$ of \eqref{eq:model} with homogeneous Neumann boundary data and $b=c_{\scriptscriptstyle N}(x,t)|\nabla v_{\scriptscriptstyle N}|^{p-1}$. Its initial datum is spatially Lipschitz and compatible with the boundary condition. Furthermore, its spatial gradient is bounded, and the source $c_{\scriptscriptstyle N}|\nabla v_{\scriptscriptstyle N}|^{p-1}$ is continuous and bounded. However, its seminorm in \eqref{eq:lip} is infinite on every right neighborhood of $-\frac{1}{2}$.
\end{prop}

\begin{remark}
From the proof of Proposition~\ref{prop:N}, we know that the initial datum is
	\begin{align*}
		\phi(x)=v_{\scriptscriptstyle N}(x,-1)=2+c_{\ast}(1-|x|^{p'})^{p'},\quad c_{\ast}=\big(16(p')^{2(p-1)}\max\{n,p'\}\big)^{-\frac{1}{p-2}},
	\end{align*}
where $p'=\frac{p}{p-1}$. Its flux is $A_{p}(\nabla\phi)=-c_{\ast}^{p-1}(p')^{2(p-1)}(1-|x|^{p'})x$, which vanishes on $\partial B_1$. Hence we have $A_{p}(\nabla \phi)\cdot\nu=0\;\,\text{on }\partial B_1.$ This is the compatibility asserted here.
\end{remark}

%\begin{remark}
%Since our results cover the case of $p=2$, it doesn't need to analyze the stability for $p$ near $2$ any more.
%
%\end{remark}

%Especially our proof holds the promise of a wide application to more types of quasilinear parabolic equations.

The remainder of the paper is organized as follows.\ Preliminary results are presented in Section \ref{SEC002}.\ Section \ref{SEC003} provides the proof of $L^2$-valued temporal H\"{o}lder regularity after any prescribed positive elapsed time.\ In Section \ref{SEC005}, we complete the proof of Theorem \ref{thm:main} by developing an elliptic reconstruction argument built upon the preceding temporal regularity. Section \ref{SEC006} is devoted to the construction of counterexamples to global higher integrability and temporal Lipschitz continuity. Finally, the appendix gives the proof of an interior elliptic estimate.

\section{Preliminary}\label{SEC002}

We adopt the following notation throughout the paper.\ Given $x_0\in\mathbb R^n$ and $\rho>0$, we denote by $B_\rho(x_0)$ the ball centered at $x_0$ with radius $\rho$. Write $B_\rho=B_\rho(0)$ and set
\begin{align*}
	Q(\rho,\tau)=B_\rho\times(-\tau,0],
	\quad \rho,\tau>0.
\end{align*}
%Accordingly, the backward parabolic cylinder with center $(x_0,t_0)$ takes the form of
%\begin{align*}
%	[(x_0,t_0)+Q(\rho,\tau)]
%	=B_\rho(x_0)\times(t_0-\tau,t_0].
%\end{align*}

A consequence of the Gagliardo–Nirenberg interpolation inequality (see e.g.\ Theorem 1.1 in \cite{DGV2012}) gives the following result.
\begin{lemma}\label{lem001}
Let $n\geq2$ and $p>1$.\ For any $B_{R}(x_{0})\subset \mathbb{R}^{n}$ and $u\in L^{2}(B_{R}(x_{0}))\cap W_{0}^{1,p}(B_{R}(x_{0}))$, we derive
\begin{align*}
\int_{B_{R}(x_{0})}|u|^{p\tilde{\chi}}dx\leq C\int_{B_{R}(x_{0})}|\nabla u|^{p}dx\bigg(\int_{B_{R}(x_{0})}|u|^{2}dx\bigg)^{\chi-1};
\end{align*}
where $C=C(n,p)$, $\chi$ and $\tilde{\chi}$ are given by
\begin{align*}
\chi=\frac{n+p}{n},\quad \tilde{\chi}=\frac{n+2}{n}.
\end{align*}

\end{lemma}

For $k\in\mathbb{R}$ and $u\in C([-1,0];L^{2}(B_{1}))\cap L^{p}((-1,0);W^{1,p}(B_{1}))$, let
\begin{align*}
	(u-k)_{+}=\max\{u-k,0\},\quad(u-k)_{-}=\max\{k-u,0\}.
\end{align*}
For an integrable function $f$ on $B_1$, write 
\begin{align*}
f_{B_1}=\dashint_{B_1}fdx=\frac{1}{|B_{1}|}\int_{B_{1}}fdx.
\end{align*} 
\begin{prop}\label{prop:truncation-poincare}
Let $u$ be the solution of \eqref{PO001} with $\Omega_{T}=Q(1,1)$ in the homogeneous Neumann case.\ Set $v=u-u_{B_1}.$\ Then there exists a constant $L_{\ast}\in(0,1)$ depending only on $n,p$ such that if $0<L\leq L_{\ast}$, for every $t\in (-1,0)$, 
	\begin{equation*}
		\|(v(\cdot,t)-k)_\pm\|_{L^p(B_1)}\leq C_1\|\nabla(v(\cdot,t)-k)_\pm\|_{L^p(B_1)},\quad \pm k\geq M_{\ast},
	\end{equation*}
	where $C_{1}=C_{1}(n,p)$ and $M_{\ast}$ is given by \eqref{Z001} below.	
\end{prop}

\begin{proof}
In the homogeneous Neumann case, testing \eqref{PO001} with $\varphi(x,t)=\zeta(t)$ for arbitrary $\zeta\in C_{0}^\infty((-1,0))$, we have $\frac{d}{dt}u_{B_1}=b_{B_1}$ in the sense of distributions. Then by the initial condition, we obtain
	\begin{equation*}
		u_{B_1}(t)=\phi_{B_1}+\int_{-1}^{t}b_{B_1}(s)ds,\quad -1\leq t\leq0.
	\end{equation*}
	Consequently, $v=u-u_{B_1}$ satisfies
	\begin{equation}\label{eq:centered-problem}
		\begin{cases}
			\partial_tv-\mathrm{div}A_{p}(\nabla v)=b-b_{B_1}(t),&\text{in }B_{1}\times(-1,0],\\
			A_{p}(\nabla v)\cdot\nu=0,
			&\text{on }\partial B_1\times (-1,0],\\
			v(\cdot,-1)=\phi-\phi_{B_1},&\text{in }B_1.
		\end{cases}
	\end{equation}
	In particular, $\nabla v=\nabla u$ and $\int_{B_1}vdx=0$. By testing \eqref{eq:centered-problem} with $v$, using a standard Steklov averaging argument, we obtain that for $t\in(-1,0),$
	\begin{align}\label{eq:energy-identity}
		\frac{1}{2}\|v(\cdot,t)\|_{L^2(B_1)}^2+\int_{-1}^{t}\!\int_{B_1}|\nabla v|^{p}dxds
		=&\frac{1}{2}\|\phi-\phi_{B_1}\|_{L^2(B_1)}^{2}+\int_{-1}^{t}\int_{B_1}bvdxds.
	\end{align}
	Combining H\"older's inequality, Poincar\'{e}'s inequality and Young's inequality, it follows that for a.e. $s\in (-1,t)$,
	\begin{align}\label{Z003}
		&\left|\int_{B_1}bvdx\right|\leq L\int_{B_{1}}|\nabla v|^{p-1}|v|dx+\int_{B_{1}}|\psi||v|dx\notag\\
		&\leq \frac{1}{8}\int_{B_{1}}|\nabla v|^{p}dx+C(p)L^{p}\int_{B_{1}}|v|^{p}dx+\|\psi\|_{L^{\frac{p}{p-1}}(B_{1})}\|v\|_{L^{p}(B_{1})}\notag\\
		&\leq\Big(\frac{1}{4}+\overline{C}(n,p)L^{p}\Big)\int_{B_{1}}|\nabla v|^{p}dx+C\int_{B_{1}}|\psi|^{\frac{p}{p-1}}dx\notag\\
		&\leq\frac{1}{2}\int_{B_{1}}|\nabla v|^{p}dx+C\int_{B_{1}}|\psi|^{\frac{p}{p-1}}dx,
	\end{align}
	provided that 
\begin{align}\label{Z005}	
0<L\leq L_{\ast}:=(4\overline{C}(n,p))^{-1/p}.	
\end{align}	
	Substituting this into \eqref{eq:energy-identity}, we deduce that for $t\in(-1,0)$,	
	\begin{align}\label{eq:energy-estimate}
		&\|v(\cdot,t)\|_{L^2(B_1)}^{2}+\int_{-1}^{t}\int_{B_1}|\nabla v|^{p}dxds\notag\\
		&\leq\|\phi-\phi_{B_1}\|_{L^2(B_1)}^{2}+C_0\int_{-1}^{t}\int_{B_{1}}|\psi|^{\frac{p}{p-1}}dxds,
	\end{align}
	where $C_{0}=C_{0}(p)$.	Let
	\begin{align}\label{Z001}
		M_{\ast}=\bigg(\frac{2}{|B_1|}\bigg)^{\frac{1}{2}}\bigg(\|\phi-\phi_{B_1}\|_{L^2(B_1)}^{2}+C_0\int_{-1}^{0}\int_{B_{1}}|\psi|^{\frac{p}{p-1}}dxds\bigg)^{\frac{1}{2}}.
	\end{align}		
	For every $k\geq M_{\ast}$, it follows from Chebyshev's inequality and \eqref{eq:energy-estimate} that for $t\in(-1,0),$
	\begin{align*}
		|\{x\in B_1:v(x,t)>k\}|\leq\frac{\|v(\cdot,t)\|_{L^2(B_1)}^2}{k^2}\leq\frac{|B_1|}{2}.
	\end{align*}
	Thus $(v-k)_+$ vanishes on a set of measure at least $|B_1|/2$. Then we have from Poincar\'{e}'s inequality that
	\begin{align*}
		\|(v(\cdot,t)-k)_+\|_{L^p(B_1)}\leq C_1\|\nabla(v(\cdot,t)-k)_+\|_{L^p(B_1)},\quad k\geq M_{\ast},
	\end{align*}
	where $C_{1}=C_{1}(n,p).$ If $k\leq-M_{\ast}$, by the same argument, we have
	\begin{align*}
		|\{x\in B_1:v(x,t)<k\}|
		\leq\frac{\|v(\cdot,t)\|_{L^2(B_1)}^2}{|k|^2}
		\leq\frac{|B_1|}{2},
	\end{align*}
	and thus,
	\begin{align*}
		\|(v(\cdot,t)-k)_-\|_{L^p(B_1)}\leq C_1\|\nabla(v(\cdot,t)-k)_-\|_{L^p(B_1)},\quad k\leq-M_{\ast}.
	\end{align*}			
	
\end{proof}

Utilizing Proposition \ref{prop:truncation-poincare}, we obtain the following global energy inequality.
\begin{lemma}\label{lem003}
Assume that $u$ solves \eqref{PO001} in $B_{1}\times(-1,0]$.\ Let $M_{\ast}$ be given in \eqref{Z001}.\ Then one can find a constant $L_{\ast}=L_{\ast}\in(0,1)$ depending only on $n,p$ such that if $0<L\leq L_{\ast}$, for $|k|\geq \max\{M_{\ast},\|\phi\|_{L^{\infty}(B_{1})}\}$ in the homogeneous Neumann case and $|k|\geq\max\{\|g\|_{L^{\infty}(\partial B_{1})},\|\phi\|_{L^{\infty}(B_{1})}\}$ in the Dirichlet case, we have
	\begin{align}\label{Z002}
		&\sup\limits_{t\in(-1,0)}\int_{B_{1}}(u-k)_{\pm}^{2}dx+\frac{1}{6}\int_{Q(1,1)}|\nabla(u-k)_{\pm}|^{p}dxdt\notag\\
		&\leq C\|\psi^{p'}\|_{L^{l_{0}}(Q(1,1))}|Q(1,1)\cap\{(u-k)_{\pm}>0\}|^{1-\frac{1}{l_{0}}},
	\end{align}
where $p'=\frac{p}{p-1}$ and $l_{0}>\frac{n+p}{p}.$	
\end{lemma}
\begin{proof}
It suffices to prove \eqref{Z002} in the homogeneous Neumann case, since the Dirichlet case can be treated by the same argument. By testing \eqref{PO001} with $\pm(u-k)_{\pm}$ and $|k|\geq \max\{M_{\ast},\|\phi\|_{L^{\infty}(B_{1})}\}$, we have from a standard Steklov averaging argument that for $t\in(-1,0)$,
\begin{align*}
&\frac{1}{2}\|(u-k)_{\pm}(\cdot,t)\|_{L^{2}(B_{1})}^{2}+\int^{t}_{-1}\int_{B_{1}}|\nabla(u-k)_{\pm}|^{p}dxds\notag\\
&=\int_{-1}^{t}\int_{B_{1}}\pm b(x,t,u,\nabla u)(u-k)_{\pm}dxds.	
\end{align*}	
In exactly the same way as in \eqref{Z003}, it follows from H\"{o}lder's inequality that 
\begin{align*}
&\left|\int_{-1}^{t}\int_{B_{1}}\pm b(x,t,u,\nabla u)(u-k)_{\pm}dxds\right|\notag\\
&\leq\frac{1}{2}\int_{-1}^{0}\int_{B_{1}}|\nabla (u-k)_{\pm}|^{p}dxds+C\int_{-1}^{0}\int_{B_{1}}|\psi|^{p'}dxds\notag\\
&\leq\frac{1}{2}\int_{-1}^{0}\int_{B_{1}}|\nabla (u-k)_{\pm}|^{p}dxds+C\|\psi^{p'}\|_{L^{l_{0}}(Q(1,1))}|Q(1,1)\cap\{(u-k)_{\pm}>0\}|^{1-\frac{1}{l_{0}}},
\end{align*}	
as $0<L\leq L_{\ast}$, where $p'=p/(p-1)$ and $L_{\ast}$ is defined by \eqref{Z005}.

\end{proof}

Making use of Lemmas \ref{lem001} and \ref{lem003}, the argument used to prove Theorem 3.2 in \cite{HLZ202601} leads to the following global $L^\infty$ estimates.
\begin{lemma}\label{LEMZ01}
Assume that $u$ is the weak solution of \eqref{PO001} with $\Omega_{T}=Q(1,1)$. Then there is a constant $L_{\ast}=L_{\ast}(n,p)\in(0,1)$ such that whenever $0<L\leq L_{\ast}$,

$(i)$ for the Dirichlet case,
\begin{align*}
	\|u\|_{L^{\infty}(Q(1,1))}\leq\max\big\{\|g\|_{L^{\infty}(\partial B_{1})},\|\phi\|_{L^{\infty}(B_{1})}\big\}+C_{2}\|\psi\|_{L^{\frac{pl_{0}}{p-1}}(Q(1,1))}^{\frac{n+p}{(p-1)(n+2)}};
\end{align*}

$(ii)$ for homogeneous Neumann case,
\begin{align*}
	\|u\|_{L^{\infty}(Q(1,1))}\leq\max\{M_{\ast},\|\phi\|_{L^{\infty}(B_{1})}\}+C_{2}\|\psi\|_{L^{\frac{pl_{0}}{p-1}}(Q(1,1))}^{\frac{n+p}{(p-1)(n+2)}},
\end{align*}
where $C_{2}=C_{2}(n,p)$, $l_{0}>\frac{n+p}{p}$, and $M_{\ast}$ is specified in \eqref{Z001}.
\end{lemma}

From Theorem 2.5 in \cite{LZ202601}, we have the following global $L^{p}_{t}W^{1,p}$ estimates. 
\begin{lemma}\label{thm03}
Assume that $u$ is the weak solution of \eqref{PO001} with $\Omega_{T}=Q(1,1)$. Then we obtain 	

$(i)$ for the Dirichlet case,
\begin{align*}
&\|u\|_{L^{p}((-1,0);W^{1,p}(B_{1}))}\notag\\
&\leq C\left(\big(\|\phi\|_{L^{\infty}(B_{1})}+\|g\|_{L^{\infty}(\partial B_{1})}\big)^{\frac{2}{p}}+\|u\|_{L^{\infty}(Q(1,1))}+\|g\|_{L^{\infty}(\partial B_{1})}\right)\notag\\
&\quad+C\Big(\|g\|_{W^{1-1/p,p}(\partial B_{1})}+\|\psi\|_{L^{\frac{p}{p-1}}(Q(1,1))}^{\frac{1}{p-1}}\Big);
\end{align*}

$(ii)$ for the homogeneous Neumann case,
\begin{align*}
&\|u\|_{L^{p}((-1,0);W^{1,p}(B_{1}))}\leq C\Big(\|\phi\|_{L^{\infty}(B_{1})}^{\frac{2}{p}}+\|u\|_{L^{\infty}(Q(1,1))}+\|\psi\|_{L^{\frac{p}{p-1}}(Q(1,1))}^{\frac{1}{p-1}}\Big),
\end{align*}
where $C=C(n,p,L)$.
	
\end{lemma}

\section{Temporal H\"{o}lder continuity}\label{SEC003}

For every
\begin{align*} 
f\in C((-1,0];L^{p}(B_{1}))\cap L^{r}((-1,0);L^{q}(B_{1})),\quad \text{with }\, r,p,q>1, 
\end{align*}
the Steklov averages are defined as follows: for $0<\varepsilon\ll1$,
\begin{align*}
	f_{\varepsilon}(x,t)=
	\begin{cases}
		\frac{1}{\varepsilon}\int^{t}_{t-\varepsilon}f(x,s)ds,& t\in(-1+\varepsilon,0),\\
		0,&t\leq-1+\varepsilon,
	\end{cases}
\end{align*}
and
\begin{align*}
	f_{\bar{\varepsilon}}(x,t)=
	\begin{cases}
		\frac{1}{\varepsilon}\int^{t+\varepsilon}_{t}f(x,s)ds,& t\in(-1,-\varepsilon),\\
		0,&t\geq-\varepsilon.
	\end{cases}
\end{align*}
For simplicity, denote
\begin{align}\label{DK02}
\begin{cases}	
\mathcal{Q}_{1}=\sup_{\tau\in(-1,0)}\|u(\cdot,\tau)\|_{L^{2}(B_{1})},\quad\mathcal{Q}_{2}=\|u\|_{L^{\infty}(Q(1,1))},\\
\mathcal{Q}_{3}=L\|\nabla u\|_{L^{p}(Q(1,1))}^{p-1}+\|\psi\|_{L^{\frac{p}{p-1}}(Q(1,1))},\quad\mathcal{Q}_{4}=\|g\|_{L^{\infty}(\partial B_{1})},
\end{cases}
\end{align}	
and
\begin{align}\label{E16}
 \mathcal{M}_{p}=
	\begin{cases}
		\mathcal{Q}_{3},&1<p<2,\\
		\mathcal{Q}_{2}^{\frac{p-2}{2(p-1)}}\mathcal{Q}_{3}^{\frac{p}{2(p-1)}},&p>2.	
	\end{cases}			
\end{align}	
We are now in a position to prove Theorem \ref{thm10}.
\begin{proof}[Proof of Theorem \ref{thm10}]
For every $0<\delta<1$ and $-1+\delta\leq s<t\leq0,$ write
\begin{align*}
h=t-s,\quad h\leq\ell \leq\frac{\delta}{2},\quad r=s-\ell,\quad \lambda=\Big(1+\frac{h}{\ell}\Big)^{\frac{1}{p-2}}.	
\end{align*}	
In light of $0<\frac{h}{\ell}\leq1,$ we obtain
\begin{align}\label{E08}
\min\big\{1,2^{\frac{1}{p-2}}\big\}\leq\lambda\leq\max\big\{1,2^{\frac{1}{p-2}}\big\},\quad |\lambda-1|\leq C_{p}\frac{h}{\ell}.
\end{align}	
Here
\begin{align*}
C_{p}=\max\big\{|p-2|^{-1},\big|2^{\frac{1}{p-2}}-1\big|\big\}\rightarrow\infty,\quad\text{as }p\rightarrow2.
\end{align*}
For $\tau\in[r,s]$, denote
\begin{align*}
T_{\lambda}(\tau)=r+\lambda^{p-2}(\tau-r),\quad\hat{u}(x,\tau)=\lambda u(x,T_{\lambda}(\tau)).
\end{align*}
Then we deduce
\begin{align}\label{E10}
T_{\lambda}(r)=r,\quad T_{\lambda}(s)=t,\quad\hat{u}(x,s)=\lambda u(x,t),\quad	\hat{u}(x,r)=\lambda u(x,r).
\end{align}	
Observe that $\hat{u}$ satisfies
\begin{equation*}
	\begin{cases}
		\partial_{\tau}\hat{u}-\mathrm{div}A_{p}(\nabla \hat{u})=\lambda^{p-1}b(x,T_{\lambda}(\tau),\lambda^{-1}\hat{u},\lambda^{-1}\nabla\hat{u}),&\text{in }B_{1}\times(r,s],\\
		\hat{u}=\lambda g(x)\;\,\text{or}\;\,
		A_{p}(\nabla \hat{u})\cdot\nu=0,&\text{on }\partial  B_{1}\times(r,s],\\
		\hat{u}(x,r)=\lambda u(x,r),&\text{in } B_{1}.
	\end{cases}
\end{equation*}
Write
\begin{align*}
F_{\lambda}(x,\tau)=\lambda^{p-1}b(x,T_{\lambda}(\tau),\lambda^{-1}\hat{u},\lambda^{-1}\nabla\hat{u})-b(x,\tau,u(x,\tau),\nabla u(x,\tau)).	
\end{align*}	
Therefore, $v(x,\tau):=\hat{u}(x,\tau)-u(x,\tau)$ solves
\begin{equation}\label{E09}
	\begin{cases}
		\partial_{\tau}v-\mathrm{div}\big(A_{p}(\nabla \hat{u})-A_{p}(\nabla u)\big)=F_{\lambda}(x,\tau),&\text{in }B_{1}\times(r,s],\\
		v=(\lambda-1)g(x)\;\,\text{or}\;\,
		(A_{p}(\nabla \hat{u})-A_{p}(\nabla u))\cdot\nu=0,&\text{on }\partial  B_{1}\times(r,s],\\
		v(x,r)=(\lambda-1) u(x,r),&\text{in } B_{1}.
	\end{cases}
\end{equation}

In the following we split the proof into two parts.

\noindent{\bf Step 1.}\ \emph{The Dirichlet case.}\ Denote $d=|\lambda-1|\mathcal{Q}_{4}$, where $\mathcal{Q}_{4}$ is defined by \eqref{DK02}. For $0<\varepsilon\ll 1$, let
\begin{align*}
 \varphi=\mathrm{sgn}(v_{\varepsilon})(|v_{\varepsilon}|-d)_{+}=\mathrm{sgn}(v_{\varepsilon})\max\{|v_{\varepsilon}|-d,0\}, 
 \end{align*}
 where $v_{\varepsilon}$ represents the Steklov average of $v$, and $\operatorname{sgn}$ denotes the sign function with $\operatorname{sgn}(i)=i/|i|$ for $i\neq0$ and $\operatorname{sgn}(0)=0$.\ Testing \eqref{E09} with the test function $\varphi_{\bar{\varepsilon}}\in W^{1,2}((r,s);L^{2}(B_{1}))\cap L^{p}((r,s);W^{1,p}_{0}(B_{1}))$, we derive that for any $\tilde{s}\in(r,s]$,
\begin{align*}
	&\int_{B_{1}}v\varphi_{\bar{\varepsilon}}dx\Big|_{r}^{\tilde{s}}+\int_{r}^{\tilde{s}}\int_{B_{1}}\big(-v\partial_{\tau}\varphi_{\bar{\varepsilon}}+(A_{p}(\nabla \hat{u})-A_{p}(\nabla u))\cdot\nabla\varphi_{\bar{\varepsilon}}\big)dxd\tau\notag\\
	&=\int_{r}^{\tilde{s}}\int_{B_{1}}F_{\lambda}(x,\tau)\varphi_{\bar{\varepsilon}}dxd\tau,
\end{align*}	
which reads that
\begin{align*}
	&\int_{r}^{\tilde{s}}\int_{B_{1}}\big(\partial_{\tau}v_{\varepsilon}\varphi+(A_{p}(\nabla \hat{u})-A_{p}(\nabla u))_{\varepsilon}\cdot\nabla\varphi\big)dxd\tau=\int_{r}^{\tilde{s}}\int_{B_{1}}(F_{\lambda}(x,\tau))_{\varepsilon}\varphi dxd\tau.
\end{align*}	
By integrating by parts and applying Lemma 3.2 in Chapter I of \cite{D1993}, we obtain
\begin{align*}
&\int_{r}^{\tilde{s}}\int_{B_{1}}\partial_{\tau}v_{\varepsilon}\varphi dxd\tau=\int_{-1}^{\tilde{s}}\int_{B_{1}}\partial_{\tau}\int^{v_{\varepsilon}}_{0}\mathrm{sgn}(\rho)(|\rho|-d)_{+}d\rho dxd\tau\notag\\
&=\frac{1}{2}\int_{B_{1}}[(|v_{\varepsilon}|-d)_{+}]^{2}(x,\tilde{s}) dx-\frac{1}{2}\int_{B_{1}}[(|v_{\varepsilon}|-d)_{+}]^{2}(x,r) dx\notag\\
&\rightarrow\frac{1}{2}\int_{B_{1}}[(|v|-d)_{+}]^{2}(x,\tilde{s}) dx-\frac{1}{2}\int_{B_{1}}[(|v|-d)_{+}]^{2}(x,r) dx,\quad\text{as }\varepsilon\rightarrow0,
\end{align*}	
and
\begin{align*}
&\int_{r}^{\tilde{s}}\int_{B_{1}}(A_{p}(\nabla \hat{u})-A_{p}(\nabla u))_{\varepsilon}\cdot\nabla\varphi dxd\tau\notag\\
&\rightarrow\int_{(B_{1}\times(r,\tilde{s}))\cap \{|v|>d\}}(A_{p}(\nabla \hat{u})-A_{p}(\nabla u))\cdot\nabla (\hat{u}-u)dxd\tau,\quad\text{as }\varepsilon\rightarrow0,
\end{align*}
and
\begin{align*}
\int_{r}^{\tilde{s}}\int_{B_{1}}(F_{\lambda}(x,\tau))_{\varepsilon}\varphi dxd\tau\rightarrow\int_{r}^{\tilde{s}}\int_{B_{1}}F_{\lambda}(x,\tau)\mathrm{sgn}(v)(|v|-d)_{+}dxd\tau,\quad\text{as }\varepsilon\rightarrow0.
\end{align*}	
From Chapter 12 in \cite{L2019}, we deduce
\begin{align*}
\int_{(B_{1}\times(r,\tilde{s}))\cap \{|v|>d\}}(A_{p}(\nabla \hat{u})-A_{p}(\nabla u))\cdot\nabla (\hat{u}-u)dxd\tau\geq0.	
\end{align*}	
Combining \eqref{E02} and \eqref{E08}, we obtain 
\begin{align*}
\|F_{\lambda}\|_{L^{\frac{p}{p-1}}(Q(1,1))}\leq C(p)\mathcal{Q}_{3},	
\end{align*}
where $\mathcal{Q}_{3}$ is specified in \eqref{DK02}. Define
\begin{align*}
\mathcal{H}:=\sup_{\tilde{s}\in(r,s)}\|(|v|-d)_{+}(\cdot,\tilde{s})\|_{L^{2}(B_{1})}.	
\end{align*}	
Consequently, using \eqref{E08}--\eqref{E10} and H\"{o}lder's inequality, it follows from these above facts that
\begin{align}\label{E11}
\mathcal{H}^{2}\leq&\|(|v|-d)_{+}(\cdot,r)\|_{L^{2}(B_{1})}^{2}+2\int_{r}^{s}\int_{B_{1}}|F_{\lambda}(x,\tau)|(|v|-d)_{+}dxd\tau\notag\\	
\leq&\mathcal{Q}^{2}_{1}|\lambda-1|^{2}+C(p)\mathcal{Q}_{3}\|(|v|-d)_{+}\|_{L^{p}(B_{1}\times(r,s))}\notag\\
\leq&\overline{C}(p)\Big(\mathcal{Q}^{2}_{1}\frac{h^{2}}{\ell^{2}}+\mathcal{Q}_{3}\|(|v|-d)_{+}\|_{L^{p}(B_{1}\times(r,s))}\Big),
\end{align}	
where $\mathcal{Q}_{1}$ and $\mathcal{Q}_{3}$ are defined in \eqref{DK02}. 

We now consider the following two subcases.

\noindent{\bf Case 1.} Let $1<p<2$. From H\"{o}lder's inequality, we obtain that for any $\tau\in(r,s)$,
\begin{align*}
\|(|v|-d)_{+}(\cdot,\tau)\|_{L^{p}(B_{1})}\leq|B_{1}|^{1/p-1/2}\|(|v|-d)_{+}(\cdot,\tau)\|_{L^{2}(B_{1})}\leq C(n,p)\mathcal{H},	
\end{align*}	
which, in combination with Young's inequality, shows that
\begin{align*}
\overline{C}(p)\mathcal{Q}_{3}\|(|v|-d)_{+}\|_{L^{p}(B_{1}\times(r,s))}\leq C\mathcal{Q}_{3}\ell^{\frac{1}{p}}\mathcal{H}\leq\frac{1}{2}\mathcal{H}^{2}+C\mathcal{Q}_{3}^{2}\ell^{\frac{2}{p}},	
\end{align*}	
where $C=C(n,p)$ and $\overline{C}(p)$ is determined by \eqref{E11}.\ Inserting this into \eqref{E11}, we derive
\begin{align}\label{E12}
\mathcal{H}\leq C(n,p)\Big(\mathcal{Q}_{1}\frac{h}{\ell}+\mathcal{Q}_{3}\ell^{\frac{1}{p}}\Big).	
\end{align}	

\noindent{\bf Case 2.} Let $p>2$. Observe first from \eqref{E08} that
\begin{align*}
(|v|-d)_{+}	\leq |v|\leq \mathcal{Q}_{2}(\lambda+1)\leq\mathcal{Q}_{2}\big(\max\big\{1,2^{\frac{1}{p-2}}\big\}+1\big),
\end{align*}	
which, together with Young's inequality, reads that
\begin{align*}
\overline{C}(p)\mathcal{Q}_{3}\|(|v|-d)_{+}\|_{L^{p}(B_{1}\times(r,s))}\leq C\mathcal{Q}_{2}^{\frac{p-2}{p}}\mathcal{Q}_{3}\ell^{\frac{1}{p}}\mathcal{H}^{\frac{2}{p}}\leq\frac{1}{2}\mathcal{H}^{2}+C\mathcal{Q}_{2}^{\frac{p-2}{p-1}}\mathcal{Q}_{3}^{\frac{p}{p-1}}\ell^{\frac{1}{p-1}}.
\end{align*}	
Substituting this into \eqref{E11}, we deduce
\begin{align*}
\mathcal{H}\leq C(p)\Big(\mathcal{Q}_{1}\frac{h}{\ell}+\mathcal{Q}_{2}^{\frac{p-2}{2(p-1)}}\mathcal{Q}_{3}^{\frac{p}{2(p-1)}}\ell^{\frac{1}{2(p-1)}}\Big).	
\end{align*}	
Combining this and \eqref{E12}, we have
\begin{align}\label{E17}
\mathcal{H}\leq C(n,p)\Big(\mathcal{Q}_{1}\frac{h}{\ell}+\mathcal{M}_{p}\ell^{\beta}\Big),	
\end{align}	
where $\beta$ and $\mathcal{M}_{p}$ are, respectively, defined by \eqref{E07} and \eqref{E16}.

Utilizing \eqref{E08}, we derive 
\begin{align}\label{E15}
\|v(\cdot,s)\|_{L^{2}(B_{1})}\leq&\|(|v|-d)_{+}(\cdot,s)\|_{L^{2}(B_{1})}+\|d\|_{L^{2}(B_{1})}\notag\\
\leq& \mathcal{H}+|\lambda-1|\mathcal{Q}_{4}|B_{1}|^{\frac{1}{2}}\leq \mathcal{H}+C(n,p)\mathcal{Q}_{4}\frac{h}{\ell}.
\end{align}	
From \eqref{E10}, we obtain that for $-1+\delta\leq s<t\leq0$,
\begin{align*}
\lambda(u(x,t)-u(x,s))=v(x,s)+(1-\lambda)u(x,s),	
\end{align*}	
which, combined with \eqref{E08} and \eqref{E15}, gives that
\begin{align}\label{E96}
\|u(\cdot,t)-u(\cdot,s)\|_{L^{2}(B_{1})}\leq&\frac{1}{\lambda}\|v(\cdot,s)\|_{L^{2}(B_{1})}+\frac{|1-\lambda|}{\lambda}\|u(\cdot,s)\|_{L^{2}(B_{1})}\notag\\
\leq&C(n,p)\Big((\mathcal{Q}_{1}+\mathcal{Q}_{4})\frac{h}{\ell}+\mathcal{M}_{p}\ell^{\beta}\Big).	
\end{align}	
Therefore, on one hand, if $0<h\leq\big(\frac{\delta}{2})^{1+\beta}$, then choosing $\ell=h^{\frac{1}{1+\beta}},$ we have
\begin{align}\label{E18}
\|u(\cdot,t)-u(\cdot,s)\|_{L^{2}(B_{1})}\leq C(n,p)\max\big\{\mathcal{Q}_{1}+\mathcal{Q}_{4},\mathcal{M}_{p}\big\}|t-s|^{\frac{\beta}{1+\beta}}.	
\end{align}	
On the other hand, if $h>\big(\frac{\delta}{2})^{1+\beta}$, we deduce
\begin{align}\label{E19}
	\|u(\cdot,t)-u(\cdot,s)\|_{L^{2}(B_{1})}\leq C(n,p)\max\big\{\mathcal{Q}_{1}+\mathcal{Q}_{4},\mathcal{M}_{p}\big\}(2/\delta)^{\beta}|t-s|^{\frac{\beta}{1+\beta}}.	
\end{align}	

\noindent{\bf Step 2.}\ \emph{The homogeneous Neumann case.}\ Write
\begin{align*}
\mathcal{K}:=\sup_{\tau\in(r,s)}\|v(\cdot,\tau)\|_{L^{2}(B_{1})}.	
\end{align*}	
Then testing \eqref{E09} with $v$ via Steklov averaging in time, it follows from \eqref{E08}--\eqref{E10} and H\"{o}lder's inequality that 
\begin{align*}
\mathcal{K}\leq&|\lambda-1|^{2}\|u(\cdot,r)\|_{L^{2}(B_{1})}^{2}+2\int_{r}^{s}\int_{B_{1}}|F_{\lambda}(x,\tau)| |v|dxd\tau\\
\leq&C\Big(\mathcal{Q}^{2}_{1}\frac{h^{2}}{\ell^{2}}+\mathcal{Q}_{3}\|v\|_{L^{p}(B_{1}\times(r,s))}\Big).
\end{align*}	
Therefore, by the same arguments as in \eqref{E17} and \eqref{E18}--\eqref{E19}, we have 
\begin{align*}
	\mathcal{K}\leq& C(n,p)\Big(\mathcal{Q}_{1}\frac{h}{\ell}+\mathcal{M}_{p}\ell^{\beta}\Big)\notag\\
\leq&
\begin{cases}
C(n,p)\max\big\{\mathcal{Q}_{1},\mathcal{M}_{p}\big\}|t-s|^{\frac{\beta}{1+\beta}},&\text{for }	0<h\leq\big(\frac{\delta}{2})^{1+\beta},\\	
C(n,p)\max\big\{\mathcal{Q}_{1},\mathcal{M}_{p}\big\}(2/\delta)^{\beta}|t-s|^{\frac{\beta}{1+\beta}},&\text{for }h>\big(\frac{\delta}{2})^{1+\beta}.	
\end{cases}			
\end{align*}	
This, in combination with \eqref{E96}, \eqref{E18}, \eqref{E19}, Lemmas \ref{LEMZ01} and \ref{thm03} shows that Theorem \ref{thm10} holds.

\end{proof}

\section{Explicit higher gradient integrability}\label{SEC005}

\subsection{Temporal defect estimate and elliptic reconstruction}

Define the nonnegative symmetric defect
\begin{equation}\label{eq:defect}
	\mathscr{D}(\xi,\zeta)=(A_{p}(\xi)-A_{p}(\zeta))\cdot(\xi-\zeta),\quad \text{for }\xi,\zeta\in\mathbb{R}^{n}.
\end{equation}
\begin{lemma}[Temporal defect bound]\label{lem:time-defect}
	For $0<\delta<1$ and $0<\tau<1-\delta$, we have
	\begin{equation}\label{eq:time-defect}
		\int_{-1+\delta}^{-\tau}\int_{B_1}
		\mathscr{D}(\nabla u(\cdot,t+\tau),\nabla u(\cdot,t))dxdt\leq C\delta^{-\beta\min\{1,2/p\}}\tau^\sigma.
	\end{equation}
	where $\beta$ and $\sigma$ are, respectively, specified in \eqref{E07} and \eqref{eq:time-exponents}.
\end{lemma}
\begin{proof}
	On $B_{1}\times(-1+\delta,-\tau)$, set
	\begin{equation*}
		\begin{aligned}
			w_\tau(x,t)&=u(x,t+\tau)-u(x,t),\\
			J_\tau(x,t)&=A_{p}(\nabla u(x,t+\tau))-A_{p}(\nabla u(x,t)),\\
			F_\tau(x,t)&=b(x,t+\tau,u(x,t+\tau),\nabla u(x,t+\tau))-b(x,t,u(x,t),\nabla u(x,t)).
		\end{aligned}
	\end{equation*}
	Then $\partial_{t}w_{\tau}-\mathrm{div} J_\tau=F_\tau$.\ Using $w_\tau$ as a test function, which is justified by a standard Steklov averaging argument, we obtain
	\begin{align}\label{F03}
		&\frac{1}{2}\|w_\tau(\cdot,-\tau)\|_{L^2(B_1)}^2
		+\int_{-1+\delta}^{-\tau}\int_{B_1}
		\mathscr{D}(\nabla u(x,t+\tau),\nabla u(x,t))dxdt\notag\\
		&=\frac{1}{2}\|w_\tau(\cdot,-1+\delta)\|_{L^2(B_1)}^2
		+\int_{-1+\delta}^{-\tau}\int_{B_1}F_\tau w_{\tau} dxdt.
	\end{align}
	By Theorem \ref{thm10}, we obtain
	\begin{align}\label{F01}
		\|w_\tau(\cdot,-1+\delta)\|_{L^2(B_1)}^2\leq C\delta^{-2\beta}\tau^{2\alpha},\quad
		\|F_\tau\|_{L^{p'}(B_1\times(-1+\delta,-\tau))}\leq2\mathcal{Q}_{3},
	\end{align}
where $\beta$ and $\mathcal{Q}_{3}$ are, respectively, given by \eqref{E07} and \eqref{DK02}. If $1<p<2$, it follows from H\"older's inequality and \eqref{F01} that
	\begin{align*}
		\|w_\tau\|_{L^p(B_1\times(-1+\delta,-\tau))}\leq C \delta^{-\beta}\tau^\alpha.
	\end{align*}
	If $p>2$, we have from Lemma \ref{LEMZ01} and \eqref{F01} that
	\begin{align*}
		\int_{-1+\delta}^{-\tau}\int_{B_1}|w_\tau|^{p}dxdt
		\leq&(2\mathcal{Q}_{2})^{p-2}\int_{-1+\delta}^{-\tau}\|w_\tau(\cdot,t)\|_{L^2(B_1)}^{2}dt\leq C\delta^{-2\beta}\tau^{2\alpha},
	\end{align*}
where $\mathcal{Q}_{2}$ is defined in \eqref{DK02}.	Thus in both cases, we derive
	\begin{equation}\label{F02}
		\|w_\tau\|_{L^p(B_1\times(-1+\delta,-\tau))}\leq C\delta^{-\beta\min\{1,2/p\}}\tau^\sigma.
	\end{equation}
	The desired estimate in \eqref{eq:time-defect} now follows by substituting \eqref{F01}--\eqref{F02} into \eqref{F03}.
\end{proof}

For an integer $N\ge1$ and $1\leq j\leq N$, define
\begin{equation}\label{eq:partition}
	h=\frac{1-\delta}{N},\quad t_j=-1+\delta+jh,\;\, 0\leq j\leq N,\quad
	I_j=(t_{j-1},t_j),
\end{equation}
and
\begin{align}\label{F05}
	\begin{cases}	
		\bar u_j(x)=h^{-1}\int_{I_j}u(x,t)dt,\quad\bar A_j(x)=h^{-1}\int_{I_j}A_{p}(\nabla u(x,t))dt,\vspace{0.3mm}\\
		\bar{b}_{j}(x)=h^{-1}\int_{I_j}b(x,t,u,\nabla u)dt,\quad d_j(x)=\frac{u(x,t_j)-u(x,t_{j-1})}{h},\\
		f_j(x)=\bar{b}_{j}(x)-d_j(x).
	\end{cases}
\end{align}
For notational convenience, denote
\begin{align}\label{AZ16}
 X=\begin{cases}
	W^{1,p}_0(B_1),&\text{for the Dirichlet case},\\
	W^{1,p}(B_1),&\text{for homogeneous Neumann case}.
\end{cases}
\end{align}
With $X$ from \eqref{AZ16}, we obtain that for any $\varphi\in X$,
\begin{equation}\label{eq:averaged-equation}
	\int_{B_1}\bar A_j\cdot\nabla\varphi dx=\int_{B_1}f_j\varphi dx.
\end{equation}
In the Neumann case, $1\in X$, so \eqref{eq:averaged-equation} gives $\int_{B_1}f_j dx=0$.\ This is the elliptic compatibility condition.

\begin{lemma}[Elliptic reconstruction]
	For each $1\leq j\leq N$, there exists $v_j\in W^{1,p}(B_1)$ such that
	\begin{equation}\label{eq:reconstruction}
		v_j-\bar u_j\in X,\quad \int_{B_1}A_{p}(\nabla v_j)\cdot\nabla\varphi dx
		=\int_{B_1}\bar A_j\cdot\nabla\varphi dx,\quad\forall\,\varphi\in X,
	\end{equation}
	where in homogeneous Neumann case we impose additionally $\int_{B_1}v_j dx=\int_{B_1}\bar u_jdx$. Moreover, we have
	\begin{equation}\label{eq:reconstruction-energy}
		\|\nabla v_j\|_{L^p(B_1)}^p\le C h^{-1}\|\nabla u\|^{p}_{L^{p}(B_{1}\times I_{j})}.
	\end{equation}
\end{lemma}
\begin{remark}\label{RE09}
	In light of \eqref{eq:averaged-equation}, we see that $v_j$ also solves $-\mathrm{div}A_{p}(\nabla v_j)=f_j$ with the selected boundary condition.
\end{remark}
\begin{proof}
	Minimize $\int_{B_1}(|\nabla v|^p/p-\bar A_j\cdot\nabla v)dx$ over $\bar u_j+W^{1,p}_0(B_1)$ in the Dirichlet case, and over $\big\{v\in W^{1,p}(B_1):\dashint_{B_1}v dx=\dashint_{B_1}\bar u_{j} dx\big\}$ in the homogeneous Neumann case. The prescribed mean fixes the additive constant in the Neumann case, and Poincar\'e's inequality provides the required control of the $W^{1,p}(B_1)$ norm. In either case, coercivity and weak lower semicontinuity yield a minimizer, which is unique by strict convexity on the admissible space.\ Its Euler equation is \eqref{eq:reconstruction}.\ Testing that equation with $v_j-\bar u_j\in X$, we have from Young's inequality and Jensen's inequality that
	\begin{align*}
		\int_{B_1}|\nabla v_j|^p dx=&\int_{B_1}A_{p}(\nabla v_j)\cdot\nabla\bar u_j dx+\int_{B_1}\bar A_j\cdot(\nabla v_j-\nabla\bar u_j)dx\\
		\leq&\frac{1}{2}\int_{B_1}|\nabla v_j|^p dx+C\int_{B_1}(|\nabla\bar u_j|^p+|\bar A_j|^{p'})dx\\
		\leq&\frac{1}{2}\int_{B_1}|\nabla v_j|^p dx+Ch^{-1}\int_{I_j}\int_{B_1}|\nabla u|^{p}dxdt,
	\end{align*}
	which implies that \eqref{eq:reconstruction-energy} holds. 
\end{proof}

Define the piecewise stationary reconstruction
\begin{equation}\label{eq:Vh}
	V_h(x,t)=v_j(x),\quad\text{for } (x,t)\in B_{1}\times I_j,\;\,1\leq j\leq N.
\end{equation}
The next identity is the central link between temporal regularity and elliptic reconstruction.

\begin{lemma}[Exact averaged defect identity]
For $1\leq j\leq N$, we have
\begin{align}\label{eq:covariance}
		&\int_{I_j}\int_{B_1}\mathscr{D}(\nabla u(x,t),\nabla v_j)dxdt\notag\\
		&=\frac{1}{2h}\int_{I_j}\int_{I_j}\int_{B_1}\mathscr{D}(\nabla u(x,t),\nabla u(x,s))dxdsdt,
\end{align}
and thus,
	\begin{equation}\label{eq:total-defect}
		\int_{B_{1}\times(-1+\delta,0)}\mathscr{D}(\nabla u,\nabla V_h)dxdt\leq C\delta^{-\beta\min\{1,2/p\}}h^{\sigma},
	\end{equation}
where $\beta$, $\sigma$ and $\mathscr{D}$ are defined in \eqref{E07}, \eqref{eq:time-exponents} and \eqref{eq:defect}, respectively.
	
\end{lemma}
\begin{proof}
	Testing \eqref{eq:reconstruction} with $v_j-\bar u_j$, we deduce that for $1\leq j\leq N,$
	\begin{equation*}
		\int_{B_1}\big(A_{p}(\nabla v_j)-\bar A_j\big)\cdot(\nabla v_j-\nabla\bar{u}_{j})dx=0,
	\end{equation*}
	which gives that
	\begin{align}\label{eq:covariance-expansion}
		&\int_{I_j}\int_{B_1}\mathscr{D}(\nabla u,\nabla v_{j})dxdt\notag\\
		&=\int_{I_j}\int_{B_1}A_{p}(\nabla u)\cdot \nabla u dxdt
		-h\int_{B_1}[\bar A_j\cdot \nabla v_{j}
		+A_{p}(\nabla v_{j})\cdot\nabla\bar{u}_{j}-A_{p}(\nabla v_{j})\cdot \nabla v_{j}]dx\notag\\
		&=\int_{I_j}\int_{B_1}A_{p}(\nabla u)\cdot \nabla u\,dx\,dt
		-h\int_{B_1}\bar A_j\cdot\nabla\bar{u}_{j}dx.
	\end{align}
	Note that
	\begin{align*}
		\int_{I_j}\int_{I_j}\mathscr{D}(\nabla u(x,t),\nabla u(x,s))dsdt=2h\int_{I_j}A_{p}(\nabla u(x,t))\cdot \nabla u(x,t)dt-2h^2\bar A_j\cdot\nabla\bar{u}_{j}.
	\end{align*}
	This, together with \eqref{eq:covariance-expansion}, shows that \eqref{eq:covariance} holds.
	
	For simplicity, let $\Phi(t,s)=\int_{B_1}\mathscr{D}(\nabla u(x,t),\nabla u(x,s))dx$.\ It then follows from Lemma \ref{lem:time-defect} and \eqref{eq:covariance} that 
	\begin{align*}
		&\int_{B_{1}\times(-1+\delta,0)}\mathscr{D}(\nabla u,\nabla V_h)dxdt=\frac{1}{h}\int_0^h\sum_{j=1}^{N}\int_{t_{j-1}}^{t_j-\tau}\Phi(s+\tau,s)dsd\tau\\
		&\leq\frac{1}{h}\int_0^h\int_{-1+\delta}^{-\tau}\Phi(s+\tau,s)dsd\tau\leq C\delta^{-\beta\min\{1,2/p\}}h^{-1}\int_0^h\tau^\sigma d\tau\leq C\delta^{-\beta\min\{1,2/p\}}h^{\sigma},
	\end{align*}
	where $\beta$ and $\sigma$ are specified in \eqref{E07}--\eqref{eq:time-exponents}. The proof is complete.
\end{proof}

%{Higher spatial norms and the gradient tail estimate} \label{sec:transfer}

\subsection{Global elliptic estimate and scalar-source lifting}

Recall that the classical Poisson kernel is given by
\begin{align*}
	P(x,\xi)=\frac{1-|x|^{2}}{|\mathbb{S}^{n-1}||x-\xi|^{n}},
	\quad x\in B_{1},\;\,\xi\in\partial B_{1}.
\end{align*}
For $g\in L^{\infty}(\partial B_{1})\cap W^{1-1/q,q}(\partial B_{1})$ with $q$  as in \eqref{eq:spatial-target}, let
\begin{align*}
	G_{0}(x)=\int_{\partial B_{1}}P(x,\xi)g(\xi)dS.
\end{align*}
Then we obtain that $G_{0}\in C^{\infty}(B_{1})$, $\Delta G_{0}=0$ in $B_{1}$, and
\begin{align}\label{POISSON02}
	\begin{cases}	
		\|G_{0}\|_{L^{\infty}(B_{1})}\leq\|g\|_{L^{\infty}(\partial B_{1})},\quad \mathrm{Tr}G_{0}=g,\\
		\|G_{0}\|_{W^{1,q}(B_{1})}\leq C(n,q)\|g\|_{W^{1-1/q,q}(\partial B_{1})},
	\end{cases}	
\end{align}
where $\mathrm{Tr}:W^{1,q}(B_{1})\longrightarrow W^{1-1/q,q}(\partial B_{1})$ is the trace operator. Here $W^{1-1/q,q}(\partial B_1)$ is the fractional Sobolev trace space. Denote
\begin{align}\label{POISSON01}
G=
\begin{cases}
G_{0},&\text{for the Dirichlet case,}\\
0,&\text{for homogeneous Neumann case.}	
\end{cases}	
\end{align}

We now give the global elliptic Calder\'on--Zygmund estimates with the same boundary condition as the reconstruction.
\begin{lemma}\label{lem:elliptic-CZ}
Let $p>1$, $p<q<\infty$, and $J\in L^{q/(p-1)}(B_1;\mathbb R^n)$.\ Suppose that $v\in W^{1,p}(B_1)$ satisfies
	\begin{equation}\label{eq:elliptic-global-weak}
		\int_{B_1}A_p(\nabla v)\cdot\nabla\varphi dx
		=\int_{B_1}J\cdot\nabla\varphi dx,
		\quad\text{for any }\varphi\in X,
	\end{equation}
where $X$ is the boundary-dependent space in \eqref{AZ16}. In the Dirichlet case we assume additionally that $v-G\in W^{1,p}_0(B_1)$ with $G$ as in \eqref{POISSON01}. Then we have
\begin{equation}\label{eq:elliptic-CZ}
		\|\nabla v\|_{L^q(B_1)}\leq C(n,p,q)\Big(\|\nabla G\|_{L^{q}(B_{1})}+\|J\|_{L^{q/(p-1)}(B_1)}^{1/(p-1)}\Big).
\end{equation}
\end{lemma}

\begin{proof}
	We first treat zero Dirichlet data and the homogeneous Neumann case, so that $v\in X$. Set $m=q/(p-1)>p'=p/(p-1)$, and define
	the closed linear subspace
	\begin{align*}
		\mathcal Z_X=\bigg\{Z\in L^{p'}(B_1;\mathbb R^n):\int_{B_1}Z\cdot\nabla\varphi dx=0\text{ for every }\varphi\in X\bigg\}.
	\end{align*}
	For $X=W^{1,p}_0(B_1)$, $\mathcal Z_X$ consists of divergence-free vector fields. For $X=W^{1,p}(B_1)$, it additionally requires zero normal flux in the weak sense. In either case, \eqref{eq:elliptic-global-weak} implies $J-A_p(\nabla v)\in\mathcal Z_X$. Let $\Pi_XJ$ be the unique minimizer of
	\begin{align*}
		Y\longmapsto\frac1{p'}\int_{B_1}|J-Y|^{p'}dx,
		\quad Y\in\mathcal Z_X.
	\end{align*}
	Using the nonlinear projection estimate in Theorem 7.2 of \cite{ISS1999} with projection exponent $p'$ and target exponent $m$, we obtain
	\begin{align}\label{F09}
		\|\Pi_XJ\|_{L^m(B_1)}\leq C(n,p,q)\|J\|_{L^m(B_1)}.
	\end{align}
	Since the field $Z:=J-A_p(\nabla v)$ belongs to $\mathcal Z_X$, we deduce that for every $Y\in\mathcal Z_X$,
	\begin{align*}
		\int_{B_1}A_{p'}(J-Z)\cdot Y dx=\int_{B_1}\nabla v\cdot Y dx=0.
	\end{align*}
	This is the Euler--Lagrange equation for the preceding minimization problem. Strict convexity therefore gives $Z=\Pi_XJ$. Consequently, we have from \eqref{F09} that
	\begin{align}\label{F10}
		\|\nabla v\|_{L^q(B_1)}^{p-1}=\|A_p(\nabla v)\|_{L^m(B_1)}=\|J-\Pi_XJ\|_{L^m(B_1)}\leq C\|J\|_{L^m(B_1)},
	\end{align}
	which proves $\|\nabla v\|_{L^q(B_1)}\leq C(n,p,q)\|J\|_{L^{q/(p-1)}(B_1)}^{1/(p-1)}$.
	
We next reduce nonzero Dirichlet case to the zero-boundary estimate. Set
	\begin{align}\label{A86}
		w=v-G\in W^{1,p}_0(B_1),\quad K=J+A_p(\nabla w)-A_p(\nabla w+\nabla G).
	\end{align}
	Then $-\mathrm{div}A_p(\nabla w)=-\mathrm{div}K$ weakly in $B_1$. Write
	\begin{align*}
		X_q=\|\nabla w\|_{L^q(B_1)},\quad
		Y_q=\|J\|_{L^m(B_1)},\quad
		Z_q=\|\nabla G\|_{L^q(B_1)},
	\end{align*}
	where $m=\frac{q}{p-1}$. On one hand, for $1<p\leq2$, using the inequality of
	\begin{align*}
		|A_p(\xi+\zeta)-A_p(\xi)|\leq C(p)|\zeta|^{p-1},\quad\text{for }\xi,\zeta\in\mathbb{R}^{n},
	\end{align*}
	we derive
	\begin{align*}
		|K|\leq |J|+C(p)|\nabla G|^{p-1}.
	\end{align*}
	This, combined with \eqref{F10}, leads to that
	\begin{align*}
		X_q^{p-1}\leq C\|K\|_{L^m(B_1)}\leq C\bigl(Y_q+Z_q^{p-1}\bigr).
	\end{align*}
	
	On the other hand, for $p>2$, we first assume that $G\in C^\infty(\overline{B_1})$. We establish $X_q<\infty$ by a finite iteration of the zero-boundary estimate. Due to the fact that
	\begin{align*}
		|A_p(\xi+\zeta)-A_p(\xi)|\leq C(p)\bigl(|\xi|^{p-2}|\zeta|+|\zeta|^{p-1}\bigr),\quad\text{for }\xi,\zeta\in\mathbb{R}^{n},
	\end{align*}
	we deduce
	\begin{align*}
		|K|\le |J|+C(p)\bigl(
		|\nabla w|^{p-2}|\nabla G|
		+|\nabla G|^{p-1}\bigr),
	\end{align*}
where $K$ is defined by \eqref{A86}. If $\nabla w\in L^r(B_1)$ for some $r\ge p$, it then follows from $\nabla G\in L^\infty(B_1)$ that
	\begin{align*}
		K\in L^s(B_1;\mathbb R^n),\quad
		s=\min\Big\{m,\frac{r}{p-2}\Big\}>p'.
	\end{align*}
	Applying \eqref{F10} with target exponent $(p-1)s>p$,
	we obtain $\nabla w\in L^{(p-1)s}(B_1)$.
	Starting from the energy exponent $r_0=p$ and iterating
	\begin{align*}
		r_{j+1}=\min\Big\{q,\frac{p-1}{p-2}r_j\Big\},\quad j\geq0,
	\end{align*}
	we reach $r_j=q$ after finitely many steps.\ Hence we have $X_q<\infty$. This preliminary argument establishes only finiteness.\ We next establish an estimate independent of $\|\nabla G\|_{L^\infty(B_1)}$. By H\"older's inequality with $1/m=(p-2)/q+1/q$, we obtain
	\begin{align*}
		\|K\|_{L^m(B_1)}&\leq Y_q+C(p)\bigl(X_q^{p-2}Z_q+Z_q^{p-1}\bigr),
	\end{align*}
	where we also used the fact of $X_q<\infty$. Consequently, a combination of \eqref{F10} and Young's inequality gives that
	\begin{align*}
		X_q^{p-1}&\leq C\bigl(Y_q+X_q^{p-2}Z_q+Z_q^{p-1}\bigr)\leq \frac12X_q^{p-1}+C\bigl(Y_q+Z_q^{p-1}\bigr),
	\end{align*}
	where $C=C(n,p,q)$. This implies that
	\begin{align*}
		X_q\leq C\bigl(Y_q^{1/(p-1)}+Z_q\bigr).
	\end{align*}
	Thus, in both cases, we have
	\begin{align}\label{AE01}
		\|\nabla v\|_{L^q(B_1)}\leq C(n,p,q)\Big(\|J\|_{L^{q/(p-1)}(B_1)}^{1/(p-1)}+\|\nabla G\|_{L^q(B_1)}\Big).
	\end{align}
	
It remains to remove the smoothness assumption on $G$ when $p>2$.\ Choose $G_k\in C^\infty(\overline{B_1})$ such that $G_k\to G$ in $W^{1,q}(B_1)$, and let $v_k\in G_k+W^{1,p}_0(B_1)$ be the unique energy solution corresponding to the same vector field $J$, satisfying that
	\begin{align*}
		\int_{B_1}A_p(\nabla v_k)\cdot\nabla\varphi dx=\int_{B_1}J\cdot\nabla\varphi dx,\quad\forall\,\varphi\in W^{1,p}_0(B_1).
	\end{align*}
	From \eqref{AE01}, we have
	\begin{align}\label{AE02}
		\|\nabla v_k\|_{L^q(B_1)}\leq C(n,p,q)\big(\|J\|_{L^{q/(p-1)}(B_1)}^{1/(p-1)}+\|\nabla G_k\|_{L^q(B_1)}\big),
	\end{align}
	which shows that $\{\nabla v_k\}$ is uniformly bounded in $L^q(B_1)$ and hence in $L^p(B_1)$. Subtracting the weak equations for $v_k$ and $v$ and testing with $(v_k-G_k)-(v-G)\in W^{1,p}_0(B_1)$, we obtain
	\begin{align}\label{AE03}
		&\int_{B_1}\bigl(A_p(\nabla v_k)-A_p(\nabla v)\bigr)\cdot(\nabla v_k-\nabla v)dx\notag\\
		&=\int_{B_1}\bigl(A_p(\nabla v_k)-A_p(\nabla v)\bigr)\cdot\nabla(G_k-G)dx\longrightarrow0,\quad \text{as }k\rightarrow\infty.
	\end{align}
	Indeed, the flux difference is bounded in $L^{p'}(B_1)$, while $\nabla(G_k-G)\to0$ in $L^p(B_1)$. From Chapter 12 in \cite{L2019}, we see that for $p>2,$ 
	\begin{align*}
		\bigl(A_p(\xi)-A_p(\zeta)\bigr)\cdot(\xi-\zeta)
		\ge 2^{2-p}|\xi-\zeta|^p,\quad\text{for }\xi,\zeta\in\mathbb{R}^{n}.
	\end{align*}
	This, in combination with \eqref{AE03}, reads that $\nabla v_k\to\nabla v$ strongly in $L^p(B_1)$.\ Together with the uniform $L^q(B_1)$ bound, this implies that after passing to a subsequence, $\nabla v_k$ converges weakly to $\nabla v$ in $L^q(B_1)$. In view of \eqref{AE02}, it then follows from weak lower semicontinuity and $G_k\to G$ in $W^{1,q}(B_1)$ that \eqref{eq:elliptic-CZ} holds.
\end{proof}

For $f\in L^r(B_1)$ with $1<r<\infty$, define the linear operator $\mathcal P$ according to the selected boundary condition. In the Dirichlet case, let
\begin{equation}\label{eq:Poisson}
	-\Delta\mathcal P f=f\quad\text{in }B_1,\qquad
	\mathcal P f=0\quad\text{on }\partial B_1.
\end{equation}
In the homogeneous Neumann case, consider
\begin{equation}\label{eq:Poisson-Neumann}
	\begin{cases}
		-\Delta\mathcal P f=f-\dashint_{B_1}f\,dx&\text{in }B_1,\\
		\partial_\nu\mathcal P f=0&\text{on }\partial B_1,\\
		\dashint_{B_1}\mathcal P f\,dx=0.&
	\end{cases}
\end{equation}
The subtraction of the mean of $f$ ensures the compatibility condition required for solvability, while the zero-mean normalization of the solution guarantees uniqueness. By the classical global linear elliptic regularity, we know
\begin{align*}
	\|\mathcal P f\|_{W^{2,r}(B_1)}
	\leq C(n,r)\|f\|_{L^r(B_1)},
\end{align*}
which, combined with the Sobolev embedding, shows that
\begin{equation}\label{eq:Poisson-bound}
	\|\nabla\mathcal P f\|_{L^m(B_1)}\leq C(n,r,m)\|f\|_{L^r(B_1)},
\end{equation}
for $1\le m\le nr/(n-r)$ if $1<r<n$, for every finite $m\geq1$ if $r=n$, and for $1\leq m\leq\infty$ if $r>n$.

For the choice of $q$ in \eqref{eq:spatial-target}, fix
\begin{equation}\label{eq:lift-exponents}
	m=\frac q{p-1},\quad r_q=\frac{nm}{n+m}=\frac{nq}{n(p-1)+q},\quad\hat{r}_{q}=\max\{2,r_q\}.
\end{equation}
Denote
\begin{equation}\label{eq:cell-masses}
	\mathcal{B}_{j}=\int_{I_j}\int_{B_1}|b(x,t,u,\nabla u)|^{p'}dxdt,\quad 1\leq j\leq N.
\end{equation}
It then follows from \eqref{eq:Poisson-bound} and H\"older's inequality that for $1\leq j\leq N,$
\begin{equation}\label{eq:forcing-lift}
	\|\nabla\mathcal P\bar{b}_{j}\|_{L^m(B_1)}^{1/(p-1)}\leq C \|\bar{b}_{j}\|_{L^{p'}(B_1)}^{1/(p-1)}
	\leq Ch^{-1/p}\mathcal{B}_{j}^{1/p},
\end{equation}
and
\begin{equation}\label{eq:endpoint-lift}
	\|\nabla\mathcal P d_j\|_{L^m(B_1)}\leq C\|d_j\|_{L^{\hat{r}_{q}}(B_1)}.
\end{equation}
The admissible range of $q$ is determined by the lifting estimate in \eqref{eq:forcing-lift}. Indeed, applying \eqref{eq:Poisson-bound} with $r=p'$ and
$m=q/(p-1)$, it requires that
\begin{align*}
\frac{q}{p-1}\leq \frac{np'}{n-p'}\;\,\Longleftrightarrow\;\, q\leq \frac{np}{n-p'},\quad \text{when }p'<n,
\end{align*}
whereas every finite $q>p$ is admissible when $p'\geq n$. Thus the argument can apply to any
\begin{align*}
\begin{cases}
	p<q\le np/(n-p'),& p'<n,\\
	p<q<\infty,& p'\ge n,
\end{cases}
\end{align*}
with the exponents defined in terms of this choice of $q$.

Set $w_j(x)=u(x,t_j)-u(x,t_{j-1})$. Then we have $d_j=w_j/h$. From Theorem \ref{thm10} and Lemma \ref{LEMZ01}, we deduce
\begin{align}\label{F07}
	\|w_j\|_{L^{\hat{r}_{q}}(B_1)}\leq(2\mathcal{Q}_{2})^{1-2/\hat{r}_{q}}\|w_j\|_{L^{2}(B_1)}^{2/\widehat{r}_q}\leq C\delta^{-2\beta/\hat{r}_{q}}h^{2\alpha/\hat{r}_{q}},
\end{align}
where $\beta$, $\mathcal{Q}_{2}$ and $\hat{r}_{q}$ are, respectively, defined by \eqref{E07}, \eqref{DK02} and \eqref{eq:lift-exponents}. From \eqref{eq:time-exponents}, we see that $2/\hat{r}_{q}=\gamma$. Then we obtain from \eqref{eq:endpoint-lift}--\eqref{F07} that
\begin{equation}\label{eq:endpoint-cost}
	\|\nabla\mathcal P d_j\|_{L^m(B_1)}^{1/(p-1)}\leq C\delta^{-\beta\gamma/(p-1)}h^{-\vartheta},
\end{equation}
where $\beta,$ $\gamma$ and $\vartheta$ are specified in \eqref{E07}, \eqref{eq:time-exponents} and \eqref{eq:gain-exponents}, respectively. Define
\begin{equation*}
	J_j=\nabla\mathcal P\bar{b}_{j}-\nabla\mathcal P d_j,\quad 1\leq j\leq N.
\end{equation*}
Then $-\mathrm{div}J_j=f_j$, where $f_j$ is defined in \eqref{F05}. In the Neumann case, we also used the fact that $\int_{B_1}f_j dx=0$. The boundary conditions defining $\mathcal P$ also imply $J_j\cdot\nu=0$ on $\partial B_1$ in the weak sense. Combining \eqref{eq:forcing-lift}--\eqref{eq:endpoint-cost}, we obtain $J_j\in L^{q/(p-1)}(B_1;\mathbb R^n)$. Moreover, $\bar b_j\in L^{p'}(B_1)$ by Jensen's inequality and $d_j\in L^\infty(B_1)$ by Lemma \ref{LEMZ01}.\ Then $\mathcal{P}\bar{b}_{j},\mathcal{P}d_{j}\in W^{2,p'}(B_1)$, and Green's identity extends by density to all tests in $X$.
For every $\varphi\in X$, combining \eqref{eq:Poisson}--\eqref{eq:Poisson-Neumann} and Remark \ref{RE09}, we obtain
\begin{equation}\label{eq:lifted-global-weak}
	\int_{B_1}A_p(\nabla v_j)\cdot\nabla\varphi dx
	=\int_{B_1}f_j\varphi dx
	=\int_{B_1}J_j\cdot\nabla\varphi dx.
\end{equation}
In the Dirichlet case, $v_j-G\in W^{1,p}_0(B_1)$, since $v_j$, $\bar u_j$, and $G$ have the same fixed trace $g$. Applying Lemma \ref{lem:elliptic-CZ} with \eqref{eq:lifted-global-weak}, it follows from \eqref{eq:forcing-lift} and \eqref{eq:endpoint-cost} that
\begin{equation}\label{eq:cell-cost}
	\|\nabla v_j\|_{L^q(B_1)}\leq C\big(h^{-1/p}\mathcal B_j^{1/p}+\delta^{-\beta\gamma/(p-1)}h^{-\vartheta}+\|\nabla G\|_{L^{q}(B_{1})}\big),
\end{equation}
where $G$ and $\mathcal{B}_{j}$ are defined in \eqref{POISSON01} and \eqref{eq:cell-masses}, respectively. Since $q/p>1$, we have
\begin{align*}
	\sum_{j=1}^{N}\mathcal{B}_{j}^{q/p}\leq\Bigg(\sum_{j=1}^{N}\mathcal{B}_{j}\Bigg)^{q/p}\leq \mathcal{Q}_{3}^{\frac{q}{p-1}},
\end{align*}
where $\mathcal{Q}_{3}$ is defined by \eqref{DK02}. Combining this with \eqref{eq:cell-cost}, we deduce from \eqref{POISSON02}, Lemmas \ref{LEMZ01} and \ref{thm03} that
\begin{align}\label{eq:reconstruction-cost}
	&\int_{B_{1}\times(-1+\delta,0)}|\nabla V_h|^{q}dxdt=\sum_{j=1}^{N}h\|\nabla v_j\|_{L^q(B_1)}^{q}\notag\\
	&\leq Ch^{1-q/p}\sum_{j=1}^{N}\mathcal{B}_{j}^{q/p}+C\delta^{-\beta\gamma q/(p-1)}Nh^{1-q\vartheta}+CNh\|\nabla G\|_{L^{q}(B_{1})}\notag\\
	&\leq C\bigl(\delta^{-\beta\min\{1,2/p\}}h^{-\kappa/p}+\delta^{-\beta\gamma q/(p-1)}h^{-q\vartheta}+1\bigr)\leq C\delta^{-\beta\gamma q/(p-1)}h^{-\mu},
\end{align}
where $\mu$ and $V_{h}$ are, respectively, defined by \eqref{eq:gain-exponents} and \eqref{eq:Vh}, and we also utilized the fact that $q-p=\kappa$ and $Nh=1-\delta\leq1$ by \eqref{eq:spatial-target} and \eqref{eq:partition}. 

\subsection{The proof of Theorem \ref{thm:main}}
Denote
\begin{equation*}
	c_{p}=\max\{4,4^{1/(p-1)}\},\quad p\in(1,2)\cup(2,\infty).
\end{equation*}
For $\xi,\zeta\in\mathbb{R}^{n}$, if $|\xi|\ge c_{p}|\zeta|$, then we have
\begin{align}\label{eq:separated-defect}
	\mathscr{D}(\xi,\zeta)\geq|\xi|^p-|\xi|^{p-1}|\zeta|-|\zeta|^{p-1}|\xi|\geq(1-c_{p}^{-1}-c_{p}^{1-p})|\xi|^p
	\geq\tfrac12|\xi|^p.
\end{align}
We are now ready to give the proof of Theorem \ref{thm:main}.
\begin{proof}[Proof of Theorem \ref{thm:main}]
	For $\Lambda\geq0$, let $E_{\Lambda}=\{(x,t)\in B_{1}\times(-1+\delta,0):|\nabla u|>\Lambda\}$. For each partition \eqref{eq:partition}, split this set according as $|\nabla V_h|\le\Lambda/c_{p}$ or $|\nabla V_h|>\Lambda/c_{p}$. A consequence of \eqref{eq:total-defect}, \eqref{eq:reconstruction-cost}, \eqref{eq:separated-defect} and Chebyshev's inequality shows that
	\begin{align}\label{eq:tail}
		|E_{\Lambda}|
		\leq&2\Lambda^{-p}\int_{B_{1}\times(-1+\delta,0)}\mathscr{D}(\nabla u,\nabla V_h)dxdt
		+c_{p}^q\Lambda^{-q}\int_{B_{1}\times(-1+\delta,0)}|\nabla V_h|^{q}dxdt\notag\\
		\leq& C\big(\delta^{-\beta\min\{1,2/p\}}\Lambda^{-p}h^\sigma+\delta^{-\beta\gamma q/(p-1)}\Lambda^{-q}h^{-\mu}\big).
	\end{align}
Using the exponents defined in \eqref{eq:spatial-target}, \eqref{eq:time-exponents} and \eqref{eq:gain-exponents}, we introduce the balancing exponent $\theta$ and its admissible partition by
	\begin{equation*}
		\theta=\frac{\kappa}{\sigma+\mu},\quad
		N=\lfloor(1-\delta)\Lambda^\theta\rfloor,
		\quad h=(1-\delta)/N,
	\end{equation*}
	where $\lfloor \cdot\rfloor$ denotes the integer part. Choose $\Lambda_0=(1-\delta)^{-1/\theta}$. Then for any $\Lambda\geq\Lambda_{0}$, we obtain that $(1-\delta)\Lambda^\theta\geq1$ and $\Lambda^{-\theta}\leq h\leq2\Lambda^{-\theta}$. The two terms in \eqref{eq:tail} are therefore bounded by $C\delta^{-\beta\gamma q/(p-1)}\Lambda^{-p-\sigma\theta}$. As $\sigma\theta=\varepsilon_*$ by \eqref{eq:gain-exponents}, we obtain that for $\Lambda\geq\Lambda_0$,
	\begin{equation}\label{eq:weak-conclusion}
		|\{(x,t)\in B_{1}\times(-1+\delta,0):|\nabla u|>\Lambda\}|\leq C\delta^{-\beta\gamma q/(p-1)}\Lambda^{-(p+\varepsilon_*)}.
	\end{equation}
	For $0<\varepsilon<\varepsilon_*$, it follows from Layer-cake formula and \eqref{eq:weak-conclusion} that if $0<L\leq L_{\ast}$ with $L_{\ast}$ as in \eqref{Z005},
	\begin{align*}
		&\int_{B_{1}\times(-1+\delta,0)}|\nabla u|^{p+\varepsilon}dxdt=(p+\varepsilon)\int_0^\infty
		\Lambda^{p+\varepsilon-1}|E_{\Lambda}|d\Lambda\notag\\
		&\leq |B_{1}\times(-1+\delta,0)|\Lambda_0^{p+\varepsilon}+C(p+\varepsilon)\delta^{-\beta\gamma q/(p-1)}\int_{\Lambda_0}^\infty\Lambda^{\varepsilon-\varepsilon_*-1}d\Lambda\notag\\
		&\leq C(\varepsilon_{\ast}-\varepsilon)^{-1}\delta^{-\beta\gamma q/(p-1)}.
	\end{align*}
The proof is complete.
\end{proof}

\section{Counterexamples to higher integrability and time regularity}\label{SEC006}
We first prove Proposition \ref{prop:counterexample}.
\begin{proof}[Proof of Proposition \ref{prop:counterexample}]
	{\bf Step 1}.\ \emph{The initial datum and its oscillatory coefficients.}\ Set
	\begin{equation*}
		\ell=\frac1{4\sqrt n},\quad D=(-\ell,\ell)^n,\quad D'=(-\ell,\ell)^{n-1}.
	\end{equation*}
	Then $\overline D\subset B_{1/2}$. Write $x=(x_1,x')$, where $x'=(x_{2},...,x_{n})$. Choose $\chi\in C_{0}^{\infty}(D')$ such that $0\leq\chi\leq1$ and $\chi\not\equiv0$. Define
	\begin{equation}\label{eq:frequencies}
		m_k=2^{k^2},\quad\lambda_k=\frac{\pi m_k}{2\ell},\quad k\geq1,
	\end{equation}
	and 
	\begin{equation*}
		\psi_k(x)=
		\begin{cases}
			\chi(x')\sin\bigl(\lambda_k(x_1+\ell)\bigr),&x\in D,\\
			0,&x\notin D.
		\end{cases}
	\end{equation*}
	Then every $\psi_k$ is Lipschitz on $\mathbb{R}^{n}$ and supported in $\overline D$. Pick $\eta\in C_{0}^{\infty}(B_{1/2})$ with $0\le\eta\le1$ and $\eta=1$ on a neighborhood of $\overline D$. Let
	\begin{equation}\label{eq:datum}
		\phi(x)=\eta(x)+\sum_{k=1}^{\infty}2^{-k}\psi_k(x).
	\end{equation}
	The series converges uniformly because
	$\|\psi_k\|_{L^\infty(B_1)}\le1$ and $\sum_{k\ge1}2^{-k}=1$.
	Consequently, $\phi$ is a continuous function compactly supported in $B_{1/2}$ with $0\leq\phi\leq2$. Note that all $m_k$ in \eqref{eq:frequencies} are even integers. Hence we have
	\begin{align*}
		\int_{-\ell}^{\ell}\sin(\lambda_k(x_1+\ell))dx_1=0,\quad k\geq1.
	\end{align*}
	From the orthogonality of the sine functions on $(-\ell,\ell)$, we deduce
	\begin{align*}
		\int_{-\ell}^{\ell}\sin(\lambda_k(x_1+\ell))\sin(\lambda_j(x_1+\ell))dx_1=
		\begin{cases}\ell,&j=k,\\0,&j\ne k.\end{cases}
	\end{align*}
	Since $\eta=1$ on $\overline D$, it follows from termwise integration of the uniformly convergent series in \eqref{eq:datum} that
	\begin{equation}\label{eq:coefficients}
		\int_{B_1}\phi\psi_{k}dx=c_{\chi}2^{-k},\quad c_\chi:=\ell\int_{D'}\chi(x')^{2}dx'>0.
	\end{equation}
	
	\noindent{\bf Step 2.}\ \emph{The bounded energy solution.}\ Set $H=L^2(B_1)$ and define
	\begin{align*}
		\mathscr E_X(f)=
		\begin{cases}
			\displaystyle\frac1p\int_{B_1}|\nabla f|^{p}dx,
			&f\in X\cap H,\\
			+\infty,&f\in H\setminus X,
		\end{cases}
	\end{align*}
	where $X$ is defined in \eqref{AZ16}.\ This functional is proper, convex, and lower semicontinuous with respect to strong convergence in $H$. Its effective domain, namely the set where it is finite, is $D(\mathscr E_X)=X\cap H$ and is dense in $H$.\ Moreover, its $H$-subdifferential is characterized by
	\begin{align*}
		a\in\partial_H\mathscr E_X(f)
		\;\,\Longleftrightarrow\;\,
		\begin{cases}
			f\in X\cap H,\quad a\in H,\\
			\displaystyle
			\int_{B_1}a\zeta dx
			=\int_{B_1}|\nabla f|^{p-2}\nabla f\cdot\nabla\zeta dx,
			\quad\zeta\in X\cap H.
		\end{cases}
	\end{align*}
	Thus the evolution inclusion $u_t+\partial_H\mathscr E_X(u)\ni0$ is precisely the $p$-Laplace equation with the boundary condition encoded by $X$.\ By using \cite[Chapter III, Theorem 3.2, pp.\ 57--58]{B1973}, there exists a unique solution $u\in C([-1,0];H)$ for every initial datum $\phi\in\overline{D(\mathscr E_X)}^{\,H}=H$.
	
	Testing the equation with the negative part and $(u-2)_+$, we derive $0\leq u\leq2$. Testing with $u$, we obtain
	\begin{equation*}
		\frac{1}{2}\|u(\cdot,t)\|_{L^2(B_1)}^{2}+\int_{-1}^{t}\int_{B_1}|\nabla u|^{p}dxds=\frac{1}{2}\|\phi\|_{L^2(B_1)}^{2}, \quad-1\leq t\leq0.
	\end{equation*}
	These tests are justified by time regularization. 
	
	\noindent{\bf Step 3.}\ \emph{A decay estimate under higher integrability.}\ For simplicity, shift the initial time by writing
	\begin{equation*}
		v(x,s)=u(x,s-1),\quad 0\leq s\leq1,\;x\in B_{1}.
	\end{equation*}
	Fix $\varepsilon>0$ and $0<\tau\leq1$. We next proceed by contradiction. Suppose that
	\begin{equation}\label{eq:assumed-bound}
		r=p+\varepsilon,\quad \mathfrak{G}:=\|\nabla v\|_{L^r(D\times(0,\tau))}<\infty.
	\end{equation}
	With the frequencies from \eqref{eq:frequencies}, define
	\begin{equation}\label{eq:primitive}
		\Psi_k(x)=
		\begin{cases}
			\displaystyle
			\chi(x')\frac{1-\cos\bigl(\lambda_k(x_1+\ell)\bigr)}
			{\lambda_k}e_1,&x\in D,\\[1mm]
			0,&x\notin D,
		\end{cases}
	\end{equation}
	where $e_1=(1,0,\ldots,0)\in\mathbb{R}^{n}$. Observe that the zero extension in \eqref{eq:primitive} is Lipschitz and has no boundary jump. Moreover, a direct differentiation gives
	\begin{equation}\label{eq:test-bounds}
		\mathrm{div}\Psi_k=\psi_k,\quad
		\|\Psi_k\|_{L^\infty(\mathbb{R}^n)}\leq2\lambda_k^{-1},
		\quad
		\|\nabla\psi_k\|_{L^\infty(\mathbb{R}^n)}\leq C\lambda_k.
	\end{equation}
	For $k\geq1$ and $s\in(0,1)$, define 
	\begin{equation*}
		I_k(s)=\int_{B_1}v(x,s)\psi_k(x)dx.
	\end{equation*}
	For every $s\in(0,\tau)$, using \eqref{eq:test-bounds}, it follows from integration by parts and H\"{o}lder's inequality that
	\begin{align} \label{eq:spatial-pairing}
		|I_k(s)|=\left|-\int_D\nabla v(x,s)\cdot\Psi_k(x)dx\right| \leq C\lambda_k^{-1}\|\nabla v(\cdot,s)\|_{L^r(D)}.
	\end{align}
	Using the tests $\psi_k$ for $k\geq1$, we obtain the following exact identity
	\begin{equation*}
		I_k(s)-I_k(0)=-\int_0^s\int_{D} A_p(\nabla v)\cdot\nabla\psi_{k}dxdt.
	\end{equation*}
	A combination of \eqref{eq:assumed-bound}, \eqref{eq:test-bounds} and H\"older's inequality shows that for $s\in(0,\tau)$,
	\begin{align}
		|I_k(s)-I_k(0)|\leq C\lambda_{k}\int_0^s\int_{D}|\nabla v|^{p-1}dxdt\leq C\lambda_k s^{1-(p-1)/r}\mathfrak{G}^{p-1}.
		\label{eq:time-bound}
	\end{align}
	Here $1-(p-1)/r=(\varepsilon+1)/r>0$.
	
	For $0<h<\tau$, we have from \eqref{eq:spatial-pairing}, \eqref{eq:time-bound} and H\"{o}lder's inequality that
	\begin{align*}
		\frac{1}{h}\int_0^h|I_k(s)|ds\leq C\lambda_k^{-1}h^{-1/r}\mathfrak{G},
	\end{align*}
	and
	\begin{align*}
		\frac{1}{h}\int_0^h|I_k(s)-I_k(0)|ds\leq C\lambda_k h^{1-(p-1)/r}\mathfrak{G}^{p-1}.
	\end{align*}
	Consequently, we obtain
	\begin{equation}\label{eq:two-costs}
		\left|\int_{B_1}\phi\psi_{k}dx\right|\leq C\big(\lambda_k^{-1}h^{-1/r}\mathfrak{G}+\lambda_k h^{1-(p-1)/r}\mathfrak{G}^{p-1}\big).
	\end{equation}
	Choose the time scale
	\begin{equation*}
		h=h_k:=\lambda_k^{-2r/(\varepsilon+2)},\quad k\geq1.
	\end{equation*}
	Since $\lambda_k\to\infty$, we have $h_k<\tau$ for all
	sufficiently large $k$. With $r=p+\varepsilon$ from
	\eqref{eq:assumed-bound}, both powers in \eqref{eq:two-costs}
	become
	\begin{align*}
		\lambda_k^{-1}h_k^{-1/r}
		&=\lambda_k^{-1+2/(\varepsilon+2)}
		=\lambda_k^{-\varepsilon/(\varepsilon+2)},\\
		\lambda_k h_k^{1-(p-1)/r}
		&=\lambda_k^{1-2(r-p+1)/(\varepsilon+2)}
		=\lambda_k^{-\varepsilon/(\varepsilon+2)}.
	\end{align*}
	Thus higher integrability would force
	\begin{equation}\label{eq:necessary-decay}
		\left|\int_{B_1}\phi\psi_{k}dx\right|\leq C(\mathfrak{G}+\mathfrak{G}^{p-1})\lambda_k^{-\varepsilon/(\varepsilon+2)},
	\end{equation}
	for all sufficiently large $k$.
	
	\noindent{\bf Step 4.}\ \emph{Contradiction.}\ In light of \eqref{eq:frequencies} and \eqref{eq:coefficients}, it follows from \eqref{eq:necessary-decay} that
	\begin{align*}
		c_{\chi}2^{-k}\leq C(\mathfrak{G}+\mathfrak{G}^{p-1})2^{-\varepsilon k^2/(\varepsilon+2)}.
	\end{align*}
	Then we derive
	\begin{align*}
		2^{\varepsilon k^2/(\varepsilon+2)-k}\leq \frac{C}{c_\chi}(\mathfrak{G}+\mathfrak{G}^{p-1}),
	\end{align*}
which is impossible as $k\to\infty$. This contradicts \eqref{eq:assumed-bound} and thus proves 
	\begin{equation*}
	\int_{-1}^{-1+\tau}\int_{D}|\nabla u|^{p+\varepsilon}dxdt=\infty,\quad
	\text{for any }\varepsilon>0\text{ and }0<\tau\leq1.
\end{equation*}	
Since
	\begin{align*}
		\int_{-1}^{0}\int_{B_1}|u|^{p+\varepsilon}dxdt\leq 2^{p+\varepsilon}|B_1|<\infty,
	\end{align*}
	we conclude that the failure comes entirely from the spatial gradient.
\end{proof}

We next prove Propositions \ref{prop:D} and \ref{prop:N}.
\begin{proof}[Proof of Proposition \ref{prop:D}]
	\noindent{\bf Step 1.}\ \emph{Common temporal filters.}\ Fix a nonzero $\zeta\in C_{0}^\infty((1,2))$ with $0\leq\zeta\leq1$, extended by zero to $\mathbb{R}$. For every integer $j\geq1$ and auxiliary time $\theta\in\mathbb{R}$, define
	\begin{equation*}
		f_j(\theta)=\zeta\bigl(4^j(1-\theta)\bigr),\qquad
		a_j(\theta)=\int_{-\infty}^{\theta}
		e^{-4^j(\theta-\sigma)}f_j(\sigma)\,d\sigma.
	\end{equation*}
	The functions $f_j$ satisfy
	$\operatorname{supp}f_j\Subset(1-2\cdot4^{-j},1-4^{-j})$, where $\Subset$ denotes compact inclusion. These intervals are mutually disjoint. Then we have
	\begin{equation}\label{eq:filters2}
		a_j'+4^ja_j=f_j,\quad 0\leq a_j\leq4^{-j},\quad |a_j'|\leq1,
		\quad a_j(\theta)=c_{\scriptscriptstyle \zeta}4^{-j}e^{-4^j(\theta-1)},
	\end{equation}
	where $\theta\geq1$ and $c_{\scriptscriptstyle \zeta}=\int_1^2e^{-r}\zeta(r)\,dr>0$.\ Fix a cube $D=(-a,a)^n\Subset B_{1/4}$, a spatial cutoff $\chi\in C_{0}^\infty(B_{1/2})$ equal to one near $\overline D$, and $\eta\in C_{0}^\infty((0,3))$ equal to one on $[\frac{1}{2},2]$. Take both cutoffs between zero and one. Let $0<\epsilon\leq1$ be a small parameter, independent of $l$ below, whose
	value will be fixed later. Fix $J=1$. For each $l\geq1$, set
	\begin{equation*}
		K_l=\frac{\pi}{a}2^{l^2+J},\quad r_{\scriptscriptstyle l}=-\frac{1}{2}+2^{-l-2},\quad
		\epsilon_{\scriptscriptstyle l}=\epsilon l^{-\frac{1}{4}},\quad \theta_l(t)=K_l^2(t-r_{\scriptscriptstyle l}),
	\end{equation*}
	and
	\begin{equation*}
		\begin{gathered}
			P_l(x,\theta)=\sum_{j=1}^{l} a_j(\theta)\sin(K_l2^jx_2),\quad
			v_{\scriptscriptstyle D}(x,t)=2+x_1+\sum_{l=1}^{\infty}w_{\scriptscriptstyle l}(x,t),
		\end{gathered}
	\end{equation*}
	where
	\begin{align*}
		w_{\scriptscriptstyle l}(x,t)=\frac{\epsilon_{\scriptscriptstyle l}}{K_l^2}\chi(x)\eta(\theta_l(t))P_l(x,\theta_l(t)).	
	\end{align*}	
	The choice $J=1$ suffices: since $D$ is compactly contained in $B_{1/4}$, we have $a<\frac{1}{4\sqrt{n}}$, and hence
	\begin{align*}
		3K_l^{-2}
		=3\left(\frac{a}{\pi}\right)^2 2^{-2l^2-2}
		<2^{-2l^2-2}.
	\end{align*}
	Then we have $3K_{l+1}^{-2}<2^{-l-3}=r_{\scriptscriptstyle l}-r_{\scriptscriptstyle l+1}$ and $r_{\scriptscriptstyle l}+3K_l^{-2}<0$. Therefore, $J_l:=(r_{\scriptscriptstyle l},r_{\scriptscriptstyle l}+3K_l^{-2})$ are disjoint, lie in $(-\frac{1}{2},0)$, and accumulate only at $-\frac{1}{2}$. Since $\sum_j a_j\leq\frac{1}{3}$, $\sum_j2^ja_j\leq1$, and $\sum_{j=1}^{l}4^ja_j\leq l$, the supremum norms on $B_{1}\times(-1,0]$ satisfy
	\begin{equation}\label{eq:block-bounds}
		\|w_{\scriptscriptstyle l}\|_{L^{\infty}}\leq C\frac{\epsilon_{\scriptscriptstyle l}}{K_l^2},\quad
		\|\nabla w_{\scriptscriptstyle l}\|_{L^{\infty}}\leq C\frac{\epsilon_{\scriptscriptstyle l}}{K_l},\quad
		\|\nabla^2w_{\scriptscriptstyle l}\|_{L^{\infty}}+\|\partial_{t}w_{\scriptscriptstyle l}\|_{L^{\infty}}\leq C\epsilon_{\scriptscriptstyle l} l.
	\end{equation}
	By disjointness, these bounds imply $\|\sum_lw_{\scriptscriptstyle l}\|_\infty+\|\nabla\sum_lw_{\scriptscriptstyle l}\|_\infty\leq C_{D,0}\epsilon$ for a fixed $C_{D,0}>0$.\ Restrict $\epsilon\leq\min\{1,(2C_{D,0})^{-1}\}$.\ This gives that
\begin{align*}	
	\frac{1}{2}\leq v_{\scriptscriptstyle D}\leq\frac{7}{2},\quad \frac{1}{2}\leq|\nabla v_{\scriptscriptstyle D}|\leq\frac{3}{2}.
\end{align*}	
	
	\noindent{\bf Step 2.}\ \emph{Cancellation and residual.}\ Let $e_1=(1,0,\ldots,0)$. Linearizing $v\mapsto\mathrm{div}A_p(\nabla v)$ about the affine function $v(x)=2+x_1$, we obtain
	\begin{align*}
	\mathcal L_p w
	:=\mathrm{div}\bigl(DA_p(e_1)\nabla w\bigr)
	=(p-1)\partial_{x_1}^2w
	+\sum_{i=2}^{n}\partial_{x_i}^2w,
	\end{align*}
where $w$ is a smooth perturbation and $DA_p$ denotes the Jacobian matrix of $A_p$.	In the following identity, all $\theta$-dependent factors are evaluated at
	$\theta=\theta_l(t)$, and $\eta'$ denotes differentiation with respect to its
	scalar argument. In particular,
	$\partial_t[\eta(\theta_l(t))]=K_l^2\eta'(\theta_l(t))$.
	Direct differentiation gives that
	\begin{align}\label{eq:Dres}
		\partial_{t}w_{\scriptscriptstyle l}-\mathcal L_pw_{\scriptscriptstyle l}=&\epsilon_{\scriptscriptstyle l}\chi\bigg(\eta\sum_{j=1}^{l}f_j\sin(K_l2^jx_2)+\eta'P_l\bigg)\notag\\
		&-\frac{\epsilon_{\scriptscriptstyle l}}{K_l^2}\eta
		\left(P_l\mathcal L_p\chi+2\partial_{x_2}\chi\,\partial_{x_2}P_l\right).
	\end{align}
	At most one $f_j$ is nonzero, so \eqref{eq:Dres} is bounded by $C\epsilon_{\scriptscriptstyle l}$. In light of \eqref{eq:block-bounds}, smoothness of the nonlinear flux near $e_1$ implies that
	\begin{equation*}
		\left|\mathrm{div}A_{p}(e_1+\nabla w_{\scriptscriptstyle l})-\mathcal L_pw_{\scriptscriptstyle l}\right|
		\leq C|\nabla w_{\scriptscriptstyle l}||\nabla^2w_{\scriptscriptstyle l}|
		\leq C\frac{\epsilon_{\scriptscriptstyle l}^2 l}{K_l}.
	\end{equation*}
	Define $F_{\scriptscriptstyle D}$ as the residual of $v_{\scriptscriptstyle D}$ on each block and zero elsewhere. Thus, we obtain
	\begin{align*}	
		\|F_{\scriptscriptstyle D}\|_{L^\infty(B_{1}\times J_l)}\leq C_{D,1}\epsilon_{\scriptscriptstyle l}\leq C_{D,1}\epsilon l^{-\frac{1}{4}}\rightarrow0,\quad\text{as }l\to\infty,
	\end{align*}
	where $C_{D,1}>0$ is independent of $l$ and $\epsilon$. Every block is smooth and vanishes near its temporal endpoints, so $F_{\scriptscriptstyle D}$ is continuous on the closed cylinder. Moreover, for any $1\leq r<\infty$, we have
	\begin{equation}\label{eq:sum-derivatives}
		\sum_{l=1}^{\infty}\int_{-1}^0
		\left(\|\partial_{t}w_{\scriptscriptstyle l}(\cdot,t)\|_{L^\infty}^{r}
		+\|\nabla^2w_{\scriptscriptstyle l}(\cdot,t)\|_{L^\infty}^{r}\right)dt
		\leq C\sum_{l=1}^{\infty}\epsilon_{\scriptscriptstyle l}^{r}l^{r}K_l^{-2}<\infty.
	\end{equation}
	Partial sums converge uniformly together with their spatial gradients.\ Moreover, \eqref{eq:sum-derivatives} permits passage in distributions with no hidden
	measure at the accumulation time. Denote
	\begin{equation*}
		c_{\scriptscriptstyle D}=\frac{F_{\scriptscriptstyle D}}
		{|\nabla v_{\scriptscriptstyle D}|^{p-1}}.
	\end{equation*}
	This coefficient is continuous and
	$\|c_{\scriptscriptstyle D}\|_{L^{\infty}}\leq2^{p-1}C_{D,1}\epsilon$.
	We now fix
	\begin{align*}
	\epsilon=\min\left\{1,\frac{1}{2C_{D,0}},
	\frac{L}{2^{p-1}C_{D,1}}\right\}>0.
	\end{align*}
	This preserves the positivity and gradient bounds and gives $\|c_{\scriptscriptstyle D}\|_\infty\leq L$.\ The choice is independent of $l$ and is kept fixed as $l\to\infty$. Thus $b=c_{\scriptscriptstyle D}|\nabla v_{\scriptscriptstyle D}|^{p-1}$ satisfies the growth condition.
	
	\noindent{\bf Step 3.}\ \emph{The endpoint fails.}\ Let
	\begin{equation}\label{eq:times}
		s_{\scriptscriptstyle l}=r_{\scriptscriptstyle l}+K_l^{-2},\quad t_{\scriptscriptstyle l}=r_{\scriptscriptstyle l}+(1+4^{-l})K_l^{-2},\quad
		d_{jl}=\frac{1-e^{-4^{j-l}}}{4^{j-l}}\geq1-e^{-1},\;\,j\leq l.
	\end{equation}
	From \eqref{eq:filters2}, we derive
	\begin{equation*}
		\frac{v_{\scriptscriptstyle D}(x,t_{\scriptscriptstyle l})-v_{\scriptscriptstyle D}(x,s_{\scriptscriptstyle l})}{t_{\scriptscriptstyle l}-s_{\scriptscriptstyle l}}
		=-\epsilon_{\scriptscriptstyle l} c_{\scriptscriptstyle \zeta}\chi(x)\sum_{j=1}^{l} d_{jl}\sin(K_l2^jx_2).
	\end{equation*}
	Since $K_la/\pi$ is an integer, it follows that $\int_{-a}^{a}\sin(K_l2^jy)\sin(K_l2^ky)\,dy=a\delta_{jk}$.\ Using $\chi=1$ on $D$ and setting $\kappa_{\scriptscriptstyle D}:=c_{\scriptscriptstyle \zeta}(1-e^{-1})[a(2a)^{n-1}]^{\frac{1}{2}}>0$, we obtain
	\begin{equation*}
	\frac{\left\|v_{\scriptscriptstyle D}(\cdot,t_{\scriptscriptstyle l})
			-v_{\scriptscriptstyle D}(\cdot,s_{\scriptscriptstyle l})\right\|_{L^{2}(B_{1})}}{|t_{\scriptscriptstyle l}-s_{\scriptscriptstyle l}|}
		\geq\kappa_{\scriptscriptstyle D}\epsilon_{\scriptscriptstyle l}\sqrt{l}
		=\kappa_{\scriptscriptstyle D}\epsilon l^{\frac{1}{4}}\longrightarrow\infty,\quad\text{as }l\to\infty.
	\end{equation*}
	
\end{proof}

\begin{proof}[Proof of Proposition \ref{prop:N}]
	\noindent{\bf Step 1.}\ \emph{A nonconstant zero-flux background.}\ Let $\Xi(x)=(1-|x|^{p'})^{p'}$ and $\overline{C}_0=(p')^{2(p-1)}\max\{n,p'\}$, where $p'=p/(p-1)$.\ Its flux is $-(p')^{2(p-1)}x(1-|x|^{p'})$, hence vanishes on $\partial B_{1}$, and
	\begin{equation*}
		\mathrm{div}\bigl(|\nabla\Xi|^{p-2}\nabla\Xi\bigr)
		=-(p')^{2(p-1)}\bigl[n-(n+p')|x|^{p'}\bigr].
	\end{equation*}
Let $V$ be the unforced homogeneous Neumann solution with initial datum $V_0(x)=2+c_{\ast}\Xi(x)$, where $c_{\ast}=(16\overline{C}_{0})^{-\frac{1}{p-2}}$. Its initial residual has norm $K_0=\overline{C}_{0}c_{\ast}^{p-1}=\frac{c_{\ast}}{16}$. By comparison with $V_0\pm K_0(t+1)$ and time translates, we obtain
	\begin{align}\label{eq:background-time}
		\begin{cases}	
			2\leq V\leq2+c_{\ast},\quad
			\|V(\cdot,t)-V(\cdot,s)\|_{L^\infty(B_{1})}\leq K_0|t-s|,\\
			\|V(\cdot,t)-V_0(\cdot)\|_{L^\infty(B_{1})}\leq K_0(t+1).
		\end{cases}	
	\end{align}
Thus $\|\partial_{t}V\|_{L^\infty(B_{1}\times(-1,0])}\leq K_0$.\ Uniqueness makes $V=V(R,t)$ radial, where $R=|x|$. From radial integration of the weak equation, we have
	\begin{equation}\label{eq:radial-flux}
		R^{n-1}|\partial_{R}V|^{p-2}\partial_{R}V
		=\int_0^R s^{n-1}\partial_{t}V(s,t)\,ds,\qquad
		|\partial_{R}V|^{p-1}\leq\frac{K_0R}{n}.
	\end{equation}
Interior gradient continuity \cite[Corollary 2.1]{L2006} and \eqref{eq:background-time} imply that $V(\cdot,t_0)$, with $t_0=-\frac{1}{2}$, is nonconstant and has a noncritical interior annulus. Choose
	\begin{align}\label{eq:noncritical-domain}
		\begin{cases}	
			\mathcal A=\{x:R_-<|x|<R_+\},\quad
			0<R_-<R_+<1,\\
			I=(t_0-\sigma,t_0+\sigma)\Subset(-1,0),
		\end{cases}
	\end{align}
	so that $|\partial_{R}V|\geq\mu_0>0$ on a neighborhood of
	$[R_-,R_+]\times\overline I$. On this noncritical cylinder, it follows from \eqref{eq:radial-flux} that $\partial_{R}^{2}V\in L^\infty$.
	The radial equation
	\begin{align*}
	\partial_{t}V-(p-1)|\partial_{R}V|^{p-2}\partial_{R}^{2}V
	=\frac{n-1}{R}|\partial_{R}V|^{p-2}\partial_{R}V
	\end{align*}
	is uniformly parabolic with H\"older coefficients there.\ Interior Schauder estimates and difference quotients then give smoothness locally on the
	same noncritical region.
	
	\noindent{\bf Step 2.}\ \emph{A local phase and globally defined perturbations.}\ On a neighborhood
	of $[R_-,R_+]\times\overline I$, define
	\begin{align}\label{eq:phase}
		\begin{cases}	
			A(R,t)=(p-1)|\partial_{R}V|^{p-2},\quad
			\mathcal L_tz=R^{1-n}\partial_{R}(R^{n-1}A\partial_{R}z),\\
			\varphi(R,t)=\int_{R_*}^{R}A(s,t)^{-\frac{1}{2}}ds,
		\end{cases}	
	\end{align}
	where $R_*\in(R_-,R_+)$. Here $\mathcal L_t$ is the spatial operator obtained by linearizing the radial diffusion about $V(\cdot,t)$. The subscript $t$ indicates the time dependence of its coefficients, not differentiation in time.\ All needed derivatives are bounded, and $A(\partial_{R}\varphi)^2=1$.\ Choose $\chi_{\scriptscriptstyle 0}\in C_{0}^\infty((R_-,R_+))$ equal to one near $R_*$. Since $\partial_{R}\varphi>0$ and $\varphi(R_*,t)=0$, compactness permits $a>0$ sufficiently small such that for every $t\in\overline I$,
	\begin{equation}\label{eq:phase-window}
		[-a,a]\Subset\varphi((R_-,R_+),t),\quad
		\{R\in(R_-,R_+):|\varphi(R,t)|\leq a\}\Subset\{\chi_{\scriptscriptstyle 0}=1\}^{\circ},
	\end{equation}
	where $\{\chi_{\scriptscriptstyle 0}=1\}^\circ$ denotes the interior of the set $\{\chi_{\scriptscriptstyle 0}=1\}$. Indeed, fix $R_-<\rho_-<R_*<\rho_+<R_+$ with $\chi_{\scriptscriptstyle 0}=1$ near $[\rho_-,\rho_+]$.\ Take $a>0$ smaller than both $-\varphi(\rho_-,t)$ and $\varphi(\rho_+,t)$ for every $t\in\overline I$, where their common lower bound is positive by compactness. Use the filters and $\eta$ from Section 2. For this separate construction, take a new small parameter $0<\epsilon\leq1$, to be fixed below independently of $l$, and choose
	\begin{align*}
		J=\max\left\{1,
		\left[\log_2\left(\frac{2\sqrt{3}\,a}{\pi\sqrt{\sigma}}\right)
		\right]\right\},
	\end{align*}
	where $\sigma$ is introduced in \eqref{eq:noncritical-domain} and $[r]$ denotes the smallest integer greater than or equal to $r$. For $l\geq1$, redefine
	\begin{align*}
		K_l=\frac{\pi}{a}2^{l^2+J},\quad
		r_{\scriptscriptstyle l}=t_0+\sigma2^{-l-2},\quad
		\epsilon_{\scriptscriptstyle l}=\epsilon l^{-\frac{1}{4}},\quad
		\theta_l(t)=K_l^2(t-r_{\scriptscriptstyle l}).
	\end{align*}
	This fixed integer satisfies $3(a/\pi)^2 2^{-2J}\leq\frac{\sigma}{4}$.
	Hence $3K_l^{-2}\leq\sigma2^{-2l^2-2}$,
	$3K_{l+1}^{-2}<\sigma2^{-l-3}=r_{\scriptscriptstyle l}-r_{\scriptscriptstyle l+1}$,
	and $r_{\scriptscriptstyle l}+3K_l^{-2}\leq t_0+\frac{3\sigma}{16}<t_0+\frac{\sigma}{4}$.
	Thus the blocks $J_l=(r_{\scriptscriptstyle l},r_{\scriptscriptstyle l}+3K_l^{-2})$
	are disjoint and contained in $(t_0,t_0+\frac{\sigma}{4})$.
	For $(R,t,\theta)\in(R_-,R_+)\times I\times\mathbb{R}$, set
	\begin{equation*}
		P_l(R,t,\theta)=\sum_{j=1}^{l}a_j(\theta)
		\sin(K_l2^j\varphi(R,t)).
	\end{equation*}
	With $R=|x|$, define 
	\begin{equation*}
		w_{\scriptscriptstyle l}(x,t)=\begin{cases}
			\dfrac{\epsilon_{\scriptscriptstyle l}}{K_l^2}\chi_{\scriptscriptstyle 0}(R)\eta(\theta_l(t))
			P_l(R,t,\theta_l(t)),&(x,t)\in\mathcal A\times I,\\
			0,&(x,t)\notin\mathcal A\times I,
		\end{cases}
	\end{equation*}
	and
	\begin{equation*}		
		\qquad v_{\scriptscriptstyle N}(x,t)=V(x,t)+\sum_{l=1}^{\infty}w_{\scriptscriptstyle l}(x,t).
	\end{equation*}
	Only $w_{\scriptscriptstyle l}$, not the phase $\varphi$, is extended by zero.\ The cutoffs make each extension smooth. With $K_\chi=\{x:|x|\in\operatorname{supp}\chi_{\scriptscriptstyle 0}\}\Subset\mathcal A$, the sum is supported in $K_\chi\times[t_0,t_0+\frac{\sigma}{4}]\Subset\mathcal A\times I$. The bounds \eqref{eq:block-bounds} and \eqref{eq:sum-derivatives} remain valid with constants depending on the fixed background solution $V$, its derivative bounds on $[R_-,R_+]\times\overline I$, the lower bound $\mu_0$, and the chosen cutoffs, but not on $l$ and $\epsilon$. Restrict $\epsilon$ to a sufficiently small interval, independently of $l$, so that
	$v_{\scriptscriptstyle N}>0$ and
	$|\nabla v_{\scriptscriptstyle N}|\geq\frac{\mu_0}{2}$ on $\mathcal A\times I$.
	Its final value is fixed after the residual estimate.
	
	\noindent{\bf Step 3.}\ \emph{Residual and coefficient: their exact domains.}\ On $\mathcal A\times I$, write
	\begin{align*}
		B_0=\partial_{R}A+\frac{n-1}{R}A,\qquad
		D_0=A\partial_{R}^{2}\varphi+B_0\partial_{R}\varphi.
	\end{align*}
	Let $\partial_{\varphi}P_l$ denote the derivative of the sine sum with respect to its phase (with $\theta$ fixed). So we have
	$\partial_{\varphi}P_l=\sum_{j=1}^{l}a_j(\theta)K_l2^j\cos(K_l2^j\varphi)$.
	In the next identity, $f_j,\eta,\eta'$ are evaluated at $\theta_l(t)$ and
	$P_l$ at $(R,t,\theta_l(t))$, while $\varphi$ retains its $(R,t)$ dependence.
	The identity $A(\partial_{R}\varphi)^2=1$ gives
	\begin{align*}
		\partial_{t}w_{\scriptscriptstyle l}-\mathcal L_tw_{\scriptscriptstyle l}
		&=\epsilon_{\scriptscriptstyle l}\chi_{\scriptscriptstyle 0}
		\bigg(\eta\sum_{j=1}^{l}f_j\sin(K_l2^j\varphi)+\eta'P_l\bigg)\\
		&\quad+\frac{\epsilon_{\scriptscriptstyle l}}{K_l^2}\eta
		\left[\bigl(\chi_{\scriptscriptstyle 0}(\partial_{t}\varphi-D_0)
		-2A(\partial_{R}\chi_{\scriptscriptstyle 0})(\partial_{R}\varphi)\bigr)
		\partial_{\varphi}P_l-(\mathcal L_t\chi_{\scriptscriptstyle 0})P_l\right],
	\end{align*}
	where $\mathcal{L}_{t}$ is defined by \eqref{eq:phase}.\ Using $\|P_l\|_{L^\infty}\leq\frac{1}{3}$,
	$\|\partial_\varphi P_l\|_{L^\infty}\leq K_l$,
	and the disjoint supports of the functions $f_j$, we find that $|\partial_{t}w_{\scriptscriptstyle l}-\mathcal L_tw_{\scriptscriptstyle l}|\leq C\epsilon_{\scriptscriptstyle l}$. On the same cylinder, it follows from Taylor expansion that
	\begin{align*}
		&\left|\mathrm{div}\bigl(|\nabla(V+w_{\scriptscriptstyle l})|^{p-2}\nabla(V+w_{\scriptscriptstyle l})
		-|\nabla V|^{p-2}\nabla V\bigr)-\mathcal L_tw_{\scriptscriptstyle l}\right|\notag\\
		&\leq
		C\bigl(|\nabla w_{\scriptscriptstyle l}||\nabla^2w_{\scriptscriptstyle l}|+|\nabla w_{\scriptscriptstyle l}|^2\bigr)
		\leq C\frac{\epsilon_{\scriptscriptstyle l}^2 l}{K_l}.
	\end{align*}
	On each $\mathcal A\times J_l$, let
	$F_{\scriptscriptstyle N}=\partial_{t}v_{\scriptscriptstyle N}
	-\mathrm{div}_x(|\nabla v_{\scriptscriptstyle N}|^{p-2}
	\nabla v_{\scriptscriptstyle N})$ and define it as zero at all remaining
	points of $\overline{B}_{1}\times[-1,0]$. The unforced background gives
	$\|F_{\scriptscriptstyle N}\|_{L^\infty(B_{1}\times J_l)}\leq C\epsilon_{\scriptscriptstyle l}=C\epsilon l^{-\frac{1}{4}}\rightarrow0$
	as $l\to\infty$, with $\epsilon$ fixed.\ Thus the residual is continuous also at $t_0$, and its support is contained in the same compact subset of $\mathcal A\times I$ as the perturbation. Define
	\begin{equation*}
		c_{\scriptscriptstyle N}(x,t)=\begin{cases}
			\dfrac{F_{\scriptscriptstyle N}(x,t)}
			{|\nabla v_{\scriptscriptstyle N}(x,t)|^{p-1}},
			&(x,t)\in\mathcal A\times I,\\[1mm]
			0,&(x,t)\notin\mathcal A\times I.
		\end{cases}
	\end{equation*}
	The denominator is at least $(\frac{\mu_0}{2})^{p-1}$ on its stated domain and the gradient is continuous there by uniform convergence.\ The numerator vanishes in a collar of the spatial and temporal boundary of $\mathcal A\times I$.\ Hence zero extension is continuous. Its norm is at most $C\epsilon$.\ We now fix $\epsilon>0$ within the previous restrictions and with $C\epsilon\leq L$. This choice may depend on the background solution $V$, annulus, cutoffs and $L$, but not on $l$, and remains fixed as $l\to\infty$. The summability in \eqref{eq:sum-derivatives} identifies the residual distributionally, with no measure at $t_0$. Near $\partial B_{1}$ the perturbation vanishes identically, so its flux equals that of $V$. The background weak equation and the convergent partial sums therefore verify the equation in the homogeneous Neumann case. The initial datum is unchanged.
	
	\noindent{\bf Step 4.}\	\emph{Accumulation of temporal difference quotients.}\ Use \eqref{eq:times}
	with these $r_{\scriptscriptstyle l},K_l$. On the annulus, freezing the phase at $s_{\scriptscriptstyle l}$, we deduce
	\begin{equation*}
		\frac{w_{\scriptscriptstyle l}(x,t_{\scriptscriptstyle l})-w_{\scriptscriptstyle l}(x,s_{\scriptscriptstyle l})}{t_{\scriptscriptstyle l}-s_{\scriptscriptstyle l}}
		=-\epsilon_{\scriptscriptstyle l} c_{\scriptscriptstyle \zeta}\chi_{\scriptscriptstyle 0}(R)
		\sum_{j=1}^{l}d_{jl}\sin(K_l2^j\varphi(R,s_{\scriptscriptstyle l}))+E_l(x),
	\end{equation*}
	where
	\begin{align*}
		E_l(x)
		=&
		\frac{\epsilon_{\scriptscriptstyle l}\chi_{\scriptscriptstyle 0}(R)}
		{K_l^2(t_{\scriptscriptstyle l}-s_{\scriptscriptstyle l})}
		\sum_{j=1}^{l}a_j(1+4^{-l})\left(
		\sin\bigl(K_l2^j\varphi(R,t_{\scriptscriptstyle l})\bigr)
		-\sin\bigl(K_l2^j\varphi(R,s_{\scriptscriptstyle l})\bigr)
		\right).
	\end{align*}
	Since $a_j(1+4^{-l})\leq4^{-j}$, we have from the mean value theorem that
	\begin{align*}
		\|E_l\|_{L^\infty(\mathcal A)}
		\leq
		C\frac{\epsilon_{\scriptscriptstyle l}}{K_l^2}
		\sum_{j=1}^{l}4^{-j}K_l2^j
		\|\partial_t\varphi\|_{L^\infty((R_-,R_+)\times I)}
		\leq C\frac{\epsilon_{\scriptscriptstyle l}}{K_l}.
	\end{align*}	
	By \eqref{eq:phase-window}, the change of variable $y=\varphi(R,s_{\scriptscriptstyle l})$
	maps the shell $|\varphi(R,s_{\scriptscriptstyle l})|<a$ onto the full interval $(-a,a)$,
	and $\chi_{\scriptscriptstyle 0}=1$ there. Set
	\begin{align*}
	W_0:=\min_{(R,t)\in[R_-,R_+]\times\overline I}
	\frac{|\mathbb S^{n-1}|R^{n-1}}{\partial_{R}\varphi(R,t)}>0,
	\qquad
	\kappa_{\scriptscriptstyle N}:=c_{\scriptscriptstyle \zeta}(1-e^{-1})(aW_0)^{\frac{1}{2}}>0.
	\end{align*}
	The Jacobian is at least $W_0$.\ Making use of orthogonality on $(-a,a)$ and \eqref{eq:background-time}, we obtain
	\begin{equation*}
		\frac{\|v_{\scriptscriptstyle N}(\cdot,t_{\scriptscriptstyle l})-v_{\scriptscriptstyle N}(\cdot,s_{\scriptscriptstyle l})\|_{L^{2}(B_{1})}}{t_{\scriptscriptstyle l}-s_{\scriptscriptstyle l}}
		\geq\kappa_{\scriptscriptstyle N}\epsilon_{\scriptscriptstyle l}\sqrt{l}-C\frac{\epsilon_{\scriptscriptstyle l}}{K_l}
		-K_0|B_{1}|^{\frac{1}{2}}\longrightarrow\infty,
	\end{equation*}
	as $l\to\infty$ with $\epsilon>0$ fixed. The constant $\kappa_{\scriptscriptstyle N}$ is independent of $l$ and $\epsilon$.
	Both times approach $-\frac{1}{2}$ from the right.
\end{proof}
\begin{remark}[The case of $p=2$]
Both counterexamples in  Propositions \ref{prop:D} and \ref{prop:N} extend to $p=2$.\ The Dirichlet construction applies without modification, since $A_2(\xi)=\xi$, $\mathcal L_2=\Delta$, and the nonlinear remainder vanishes identically. For the Neumann construction, replace the background in Step~1 by an explicit solution of the heat equation.\ Choose a nonconstant radial Neumann eigenfunction $\Psi\in C^\infty(\overline{B}_1)$ satisfying
	\begin{align*}
	-\Delta\Psi=\lambda\Psi\quad\text{in }B_1,
	\quad
	\partial_\nu\Psi=0\quad\text{on }\partial B_1,
	\quad
	\|\Psi\|_{L^\infty(B_1)}=1,
	\end{align*}
	where $\lambda>0$. Let
	\begin{align*}
	V(x,t)=2+e^{-\lambda(t+1)}\Psi(x), \quad (x,t)\in\overline{B}_1\times[-1,0].
	\end{align*}
Then we have $V_t-\Delta V=0$, $\partial_\nu V=0$, and
	\begin{align*}
	1\leq V\leq3,
	\quad
	\|V(\cdot,t)-V(\cdot,s)\|_{L^\infty(B_1)}
	\leq\lambda|t-s|.
	\end{align*}
Its initial datum $V(\cdot,-1)=2+\Psi$ is smooth and satisfies the homogeneous Neumann condition. Since $\Psi$ is radial and nonconstant, its radial derivative is nonzero at some radius in $(0,1)$. Thus a noncritical cylinder $\mathcal A\times I$ around $t_0=-1/2$ can be chosen as in \eqref{eq:noncritical-domain}. The subsequent perturbation argument applies with this background and $K_0=\lambda$. In \eqref{eq:phase}, we now have
	\begin{align*}
	A\equiv1,
	\quad
	\varphi(R,t)=R-R_*.
	\end{align*}
Consequently, both the nonlinear remainder and the phase error $E_l$ vanish identically. The residual estimates and the construction of the continuous coefficient $c_{\scriptscriptstyle N}$ remain valid. With $\epsilon>0$ fixed sufficiently small, the same orthogonality argument yields
	\begin{align*}
	\frac{
		\|v_{\scriptscriptstyle N}(\cdot,t_l)
		-v_{\scriptscriptstyle N}(\cdot,s_l)\|_{L^2(B_1)}
	}{t_l-s_l}
	\geq
	\kappa_{\scriptscriptstyle N}\epsilon l^{1/4}
	-\lambda|B_1|^{1/2}
	\longrightarrow\infty,\quad \text{as }l\rightarrow\infty.
	\end{align*}
Hence the failure of temporal Lipschitz continuity also holds for $p=2$.
\end{remark}

\appendix
\section{An interior elliptic estimate}
For comparison, we also present the corresponding interior estimate together with an explicit covering argument. This appendix requires no boundary condition and is not used to obtain the global conclusions above.

\begin{lemma}[Interior elliptic estimate]\label{interior-elliptic-CZ}
	Let $p>1$ and $p<q<\infty$. Suppose that $v\in W^{1,p}(B_1)$ and $J\in L^{q/(p-1)}(B_1;\mathbb R^n)$ satisfy $-\mathrm{div}A_{p}(\nabla v)=-\mathrm{div} J$ weakly in $B_1$. Then for any $0<R<1$, we have
	\begin{equation}\label{eq:interior-elliptic-CZ}
		\|\nabla v\|_{L^q(B_R)}\leq C(n,p,q)\Big((1-R)^{-n(1/p-1/q)}\|\nabla v\|_{L^p(B_1)}+\|J\|_{L^{q/(p-1)}(B_1)}^{1/(p-1)}\Big).
	\end{equation}
	
\end{lemma}
\begin{remark}
	The interior estimate in \eqref{eq:interior-elliptic-CZ} can be also obtained by using the local nonlinear Calder\'{o}n--Zygmund estimate in Theorem 1.3 of \cite{BCDKS2018}, local Fefferman--Stein inequality and a standard covering argument.	
\end{remark}

\begin{proof}
	Utilizing Theorem 2.2 in \cite{M2011}, we deduce that for any $B_{2s}(x_0)\Subset B_1$,
	\begin{align}\label{eq:elliptic-two-ball}
		\bigg(\dashint_{B_s(x_0)}|\nabla v|^{q}dx\bigg)^{\frac{1}{q}}\leq C \bigg(\dashint_{B_{2s}(x_0)}|\nabla v|^{p}dx\bigg)^{\frac{1}{p}}+C\bigg(\dashint_{B_{2s}(x_0)}|J|^{\frac{q}{p-1}}dx\bigg)^{\frac{1}{q}},
	\end{align}
	where $C=C(n,p,q).$
	
	The covering with bounded overlap follows from \cite[Lemma 2.2]{CF2019} and \cite[Lemma 4.2.6 and Remark 4.2.7]{V2018}. For completeness, we give the detailed construction below. Fix $0<R<1$ and set
	\begin{equation*}
		d=1-R,\quad\rho=\frac{d}{8}.
	\end{equation*}
	Here $d$ represents the distance from $\overline{B_{R}}$ to $\partial B_1$, and $\rho$ denotes the covering radius. Choose points in $\overline{B_{R}}$ successively, requiring each new point to have distance at least $\rho$ from all previously chosen points, and stop when no further point can be added. This procedure terminates after finitely many steps. Indeed, if $\ell$ points have been selected, then the open balls $B_{\rho/2}(x_i)$, $1\le i\le\ell$, are pairwise disjoint and are contained in $B_{R+\rho/2}$. Write $\omega_n=|B_1|$. Comparison of their volumes gives that
	\begin{align*}
		\ell\omega_n\left(\frac{\rho}{2}\right)^n\leq \omega_n\left(R+\frac{\rho}{2}\right)^n,\quad\ell\leq\bigg(1+\frac{2R}{\rho}\bigg)^n.
	\end{align*}
	Denote the resulting family by
	$\{x_i\}_{i=1}^{N_R}\subset\overline{B_{R}}$, where $N_R$ is its number of points. This family is \emph{$\rho$-separated}, meaning that
	\begin{equation}\label{eq:elliptic-center-separation}
		|x_i-x_k|\geq\rho,\quad\text{whenever }i\neq k.
	\end{equation}
	It is \emph{maximal with respect to inclusion}: no point of $\overline{B_{R}}$ can be added while preserving \eqref{eq:elliptic-center-separation}. Maximality here does not mean that the family has the largest possible cardinality.
	
	If a point $y\in\overline{B_{R}}$ belongs to none of the balls $B_\rho(x_i)$, then $|y-x_i|\ge\rho$ for every $i$. So $y$ could be added to the family while preserving $\rho$-separation, contradicting its maximality. Consequently, we obtain
	\begin{equation*}
		\overline{B_{R}}
		\subset\bigcup_{i=1}^{N_R}B_\rho(x_i).
	\end{equation*}
	Combining this covering with the preceding volume comparison, it follows that
	\begin{equation*}
		\max\left\{1,\left(\frac{R}{\rho}\right)^n\right\}\leq N_R\leq\left(1+\frac{2R}{\rho}\right)^n\leq\left(\frac{1+15R}{1-R}\right)^n.
	\end{equation*}
	In fact, the lower bound follows from $|B_R|\leq N_R|B_\rho|$ and $N_R\geq1$.\ The value of $N_R$ may depend on the selected family. Note that 
	\begin{align*}
		\overline{B_{2\rho}(x_i)}\subset\overline{B_{R+2\rho}}\subset B_1,\quad R+2\rho=\frac{1+3R}{4}<1.
	\end{align*}
	In particular,
	\begin{align*}
		\operatorname{dist}\bigl(\overline{B_{2\rho}(x_i)},\partial B_1\bigr)\geq\frac{3}{4}(1-R)>0.
	\end{align*}
	Then the local elliptic estimate in \eqref{eq:elliptic-two-ball} can apply on each $B_{2\rho}(x_i)$.
	
	To estimate the overlap, fix $x\in\mathbb R^n$ and define
	\begin{align*}
		\mathcal I(x)
		=\{i\in\{1,\ldots,N_R\}:x\in B_{2\rho}(x_i)\}.
	\end{align*}
	Thus $\#\mathcal I(x)$ is the number of enlarged balls containing $x$, where $\#$ denotes cardinality.\ For every $i\in\mathcal I(x)$, $B_{\rho/2}(x_i)\subset B_{5\rho/2}(x)$.\ Since these smaller balls are pairwise disjoint, then we have
	\begin{align*}
		\#\mathcal I(x)\omega_n\left(\frac{\rho}{2}\right)^n=\sum_{i\in\mathcal I(x)}|B_{\rho/2}(x_i)|\leq \omega_n\Big(\frac{5\rho}{2}\Big)^n.
	\end{align*}
	Therefore, it follows that
	\begin{equation}\label{eq:elliptic-overlap}
		\sum_{i=1}^{N_R}\mathbf 1_{B_{2\rho}(x_i)}(x)=\#\mathcal I(x)\leq 5^n,
		\quad \text{for }x\in\mathbb R^n,
	\end{equation}
	where $\mathbf 1_E$ is the indicator function of the set $E$.\ Especially for every nonnegative $w\in L^1(B_1)$, we have from \eqref{eq:elliptic-overlap} that
	\begin{equation}\label{eq:elliptic-overlap-integral}
		\sum_{i=1}^{N_R}\int_{B_{2\rho}(x_i)}w(x)dx=\int_{B_1}w(x)\sum_{i=1}^{N_R}\mathbf 1_{B_{2\rho}(x_i)}(x)dx\leq 5^n\int_{B_1}w(x)dx.
	\end{equation}
	Applying \eqref{eq:elliptic-two-ball} with $s=\rho$, we obtain
	\begin{equation}\label{eq:elliptic-cover-local}
		\int_{B_\rho(x_i)}|\nabla v|^{q}dx\leq C\bigg(\rho^{n-nq/p}\bigg(\int_{B_{2\rho}(x_i)}|\nabla v|^pdx\bigg)^{\frac{q}{p}}+\int_{B_{2\rho}(x_i)}|J|^{\frac{q}{p-1}}dx\bigg),
	\end{equation}
	where $C=C(n,p,q)$.\ Since $q/p>1$, one has $\sum_i a_i^{q/p}\le(\sum_i a_i)^{q/p}$ for $a_{i}\geq0$. This, in combination with  \eqref{eq:elliptic-overlap-integral}--\eqref{eq:elliptic-cover-local}, leads to that
	\begin{align*}
		\int_{B_R}|\nabla v|^{q}dx&\leq C(n,p,q)\bigg(\rho^{n-nq/p}\left(\int_{B_1}|\nabla v|^{p}dx\right)^{\frac{q}{p}}+\int_{B_1}|J|^{\frac{q}{p-1}}dx\bigg),
	\end{align*}
	which reads that \eqref{eq:interior-elliptic-CZ} holds. 
\end{proof}

%
%%
%\noindent{\bf{\large Conflict of interest.}} The authors declare that they have no conflict of interest.
%\noindent{\bf{\large Data Availability Statement.}} The data used to support the findings of this study are available from the corresponding author upon request.

\noindent{\bf{\large Acknowledgements.}} Z. Zhao was partially supported by NSFC (No.\ 12501254).

%\noindent{\bf{\large Statements.}}

%\noindent{\bf{\large Acknowledgements.}}

\end{document}